\documentclass[final,leqno,letterpaper]{etna}

\usepackage{amssymb}
\usepackage{amsmath}
\usepackage{bm}
\usepackage{graphicx}
\usepackage{mathtools}
\usepackage{dsfont}
\usepackage[space,noadjust,nocompress]{cite}

\numberwithin{mythm}{section}

\usepackage{color}
\usepackage[dvipsnames]{xcolor}
\usepackage[caption=false]{subfig}
\usepackage{pgfplots}
\usepackage{comment}
\pgfplotsset{compat=1.18}

\definecolor{myblue}{rgb}{0.3, 0.3, 4.0}

\definecolor{mygreen}{rgb}{0.1, 0.5, 0.1}

\definecolor{markgreen}{rgb}{0,0.5,0}
\definecolor{markblue}{rgb}{0,0.1,0.9}
\definecolor{markred}{rgb}{0.85,0,0}

\makeatletter
\newenvironment{proofof}[1]{\par{\textit{Proof of
#1.}}\ignorespaces\global\@endproofsetfalse}{\endproof}
\makeatother

\setbibdata{1}{xx}{46}{2017} 

\hypersetup{%
  pdftitle={Galerkin methods for 3D unconventional superconductors},
  pdfauthor={Andreas A. Buchheit},
  pdfkeywords={}
}

\title{Galerkin methods for unconventional superconductors with
  long-range interactions on 3D lattices
  \thanks{%
    Received... Accepted... Published online on... Recommended by....
}}

\author{Andreas A. Buchheit\footnotemark[2]~\footnotemark[3]
  \and Torsten Ke{\ss}ler\footnotemark[4]~\footnotemark[5] \and Sergej
Rjasanow\footnotemark[2]}

\shorttitle{Galerkin methods for 3D unconventional superconductors}
\shortauthor{A.~BUCHHEIT AND T.~KESSLER AND S.~RJASANOW}

\begin{document}

\maketitle

\renewcommand{\thefootnote}{\fnsymbol{footnote}}

\footnotetext[2]{Saarland University, 66123 Saarbr\"ucken, Germany}
\footnotetext[3]{ETH Z\"urich, 8092 Z\"urich, Switzerland}
\footnotetext[4]{Eindhoven University of Technology, 5612 Eindhoven,
Netherlands}
\footnotetext[5]{Simkinetic B.V., Eindhoven,
Netherlands}

\begin{abstract}
  We develop an analysis and a B-Spline Galerkin solution for the
  matrix-valued Bardeen--Cooper--Schrieffer (BCS) equation for the
  mathematical description of unconventional superconductors  with
  general long-ranged interaction kernels on arbitrary Bravais
  lattices at finite temperature. To this end, we first define
  the BCS free energy functional on lattices, show its boundedness
  from below, and
  characterise solutions of the gap equation as stationary points of
  the functional. For non-repulsive interaction kernels, we prove
  that fixed-point iteration never increases the free energy and
  drives the residual of the gap equation to zero.  Our
  formulation separates the nonlinear BCS equation into a linear
  pseudodifferential equation within the Brillouin zone of the
  crystalline structure and a nonlinear functional term. The
  algebraically decaying tails of the interaction kernel are
  described by Epstein zeta functions, accompanied by a finite
  Fourier series for the short-ranged contribution. Since in our
  discretisation the interaction kernel only enters through its
  values at lattice sites, the method applies to any power-law
  exponent $\nu>0$, including exponents equal to or below the lattice
  dimension, i.e.~dipolar interactions in three-dimensional lattices.
  In the B-spline Galerkin basis, the arising singular convolutions
  correspond to block-circulant matrices with eigenvalues given by
  efficently computable lattice sums, applied without quadrature at
  FFT cost. Furthermore, through engineering a suitable dispersion
  for a prescribed gap, we obtain a set of non-trivial analytic
  singlet and triplet solutions to the BCS gap equation with
  power-law interactions on the integer lattice at zero temperature,
  the first ones, to the best of our knowledge.
  Exploiting the convolution periodical structure of the BCS gap
  equation, our numerical method combines periodic B-splines of
  arbitrary order, circulant and block-circulant matrices, and FFT,
  yielding log-linear work per iteration step and linear memory.
  Convergence is accelerated using Broyden-type quasi-Newton updates.
  Benchmarks against the analytical solutions and against highly
  resolved numerical data show reliable convergence, with rates
  determined by the B-spline degree and the regularity of the solution.
  Finally, we apply our method to a nodal superconductor with
  long-range interactions on the three-dimensional fcc lattice,
  obtaining a chiral state with a line node and two topologically
  non-trivial Weyl point nodes. On 16 cores, the method requires less
  than two minutes to obtain a residual of $1.5\times10^{-13}$ for
  this three-dimensional problem with $6.3$ million complex unknowns.
\end{abstract}



\section{Introduction}
\label{sec:intro}

Superconductivity is the quantum phenomenon where a material exhibits
vanishing electric resistance when cooled below a critical
temperature $T_c$. Although discovered already in mercury by Onnes in
1911, the first breakthrough in its microscopic understanding
occurred only in 1957 in the work of Bardeen, Cooper, and Schrieffer
\cite{BCS} that explained superconductivity as the result of the
formation of a coherent state of bosonic electron pairs, so called
Cooper pairs bound together by an effectively attractive interaction
mediated by lattice vibrations, called phonons. The pairing between
electrons is quantified by the superconducting gap $\Delta(\bm x)\in
\mathds C^{2\times 2}$, where $\bm x$ is the electron separation or
its Fourier transform $\hat \Delta(\bm p)$, where $\bm p$ is wave
vector, in general a matrix valued quantity, that describes the
binding energy of the Cooper pair, the spatial or momentum structure,
as well as the spin structure of the superconducting state. A
non-zero gap signals the presence of superconductivity.  To date, we
are aware that the structure of the gap in wave vector and spin
degrees of freedom can become very complicated and determining the
gap numerically is therefore a challenging task, especially if the
gap is slowly decaying in real space, and if the interaction between
electrons is long-ranged \cite{buchheit2023exact,haink2025non}.

In this work, we present an efficient and general numerical method
for efficiently solving the nonlinear BCS gap equation for electrons
on arbitrary $d$-dimensional Bravais lattices, based on B-spline
Galerkin discretisation combined with a pseudo-differential operator
treatment, including both singlet and triplet channels, finite
temperature, as well as attractive electron-electron interactions
that may include general power-law kernels accompanied by a compactly
supported near-field term. The method offers log-linear scaling in
both, work and memory, and permits realistic simulations of
unconventional and topological superconductivity, in one-, two-, and
three-dimensional lattices. In particular, it opens the door to
realistic simulations of three-dimensional nodal states.

We begin by providing a brief account of the most relevant
experimental discoveries and theoretical advances in unconventional
superconductors and proceed with discussing the relevance of
long-ranged interactions in current superconductivity research. We
then provide an overview of previous mathematical and numerical works
in the field.

Until the mid-1970s, the superconducting gap was assumed to have a
simple structure, smooth and with vanishing relative orbital angular
momentum with opposite spins of the participating electrons. This is
referred to as $s$-wave pairing, following the naming of atomic
orbital symmetries, and the resulting materials are called
conventional superconductors. The first discovery of unconventional
Cooper pairing, exhibiting a more complicated structure of the gap
function, was made by Osheroff in 1971
\cite{osheroff1972evidence} in superfluid
$\,^3$He, whose A- and B- phases are examples for odd momentum and
spin-triplet pairing \cite{leggett1975theoretical}, where also the
diagonal elements of the gap matrix can be non-zero. In solids,
unconventional superconductivity was discovered in heavy Fermion
compounds in 1979 in CeCu$_2$Si$_2$
\cite{steglich1979superconductivity}, followed by the discovery of
unconventional pairing in organic conductors
\cite{jerome1980superconductivity}. In those cases, electron-phonon
coupling cannot explain the appearance of superconductivity, as it is
dominated by the Coulomb repulsion. Here, other forms of effective
interactions play a role in creating the superconducting state, many
of which are closely related to the presence of magnetic ordering in
the vicinity of the superconducting transition. Unconventional
pairing received broad attention with the discovery of
high-temperature superconductivity in cuprate compounds in 1986
\cite{bednorz1986possible}, where the pairing interaction arises not
from phonon but from electronic correlations
\cite{scalapino2012common}. Other examples for unconventional pairing
include iron-based superconductors \cite{kamihara2008iron}, where the
three-dimensional structure is particularly important, see
\cite{stewart2011superconductivity} for a review.  For general
reviews on the phenomenology of unconventional superconductivity and
an analysis of the arising symmetries, see
\cite{sigrist1991phenomenological,stewart2017unconventional,norman2011challenge}.
An odd gap function in the wave vector argument, referred to as
odd-parity pairing, implies the existence of topologically protected
boundary modes
\cite{chiu2016classification,kallin2016chiral,kitaev2001unpaired,sato2017topological,schnyder2008classification},
topological surface states exist even in the presence of nodes
(points in momentum where the gap matrix vanishes)
\cite{schnyder2015topological}, and resulting topological Majorana
zero modes are currently sought for as building blocks of
error-resistant qubits for quantum computing purposes
\cite{alicea2012new,lutchyn2018majorana}. A current candidate for
odd-parity superconductivity is the material UTe$_2$
\cite{ran2019nearly,li2025observation}.

Long-ranged electron-electron interactions are a potential route to
achieving such exotic states of matter. Their realisation ranges from
dipolar and Rydberg-dressed gases
\cite{baranov2012condensed,cooper2009stable} to cavity-mediated
interactions \cite{PhysRevLett.125.053602}. Furthermore, effectively
long-ranged electron-electron are created when approaching a quantum
critical point. Here, fluctuations of the order parameter mediate an
effective interaction whose range diverges at the critical point
\cite{lohneysen2007fermi}. In previous works, we have shown that
interactions decaying as a power-law of the distance can yield
mixed-parity and topologically non-trivial superconducting states
\cite{buchheit2023exact}, where the interaction kernel is described
in terms of Epstein zeta functions \cite{epstein1903theorieI,
epstein1903theorieII, buchheit2024epstein}. We have further shown
that topological Majorana edge modes that arise from long-range
interactions in a $p$-wave superconducting string inherit the
long-range interactions even if their wave functions do not overlap
\cite{haink2025non}.

Compared to the vast body of work in condensed matter physics, research
on the mathematical properties and, in particular, the efficient
numerical simulation of unconventional superconductors are much more
scarcer. The original gap equation, a nonlinear Hammerstein-type
integral equation for the lattice Fourier transform of the gap matrix
$\hat\Delta$, was first provided by Bardeen, Cooper, and Schrieffer
in \cite{bardeen1957theory}. The mathematical structure of the
Bardeen-Cooper-Schrieffer (BCS) model was analyzed in the beginning
of the 1960s \cite{haag1962mathematical}, based on the Bogoliubov
transformation \cite{bogoliubov1958new}. Existence and uniqueness of
positive solutions of the gap equation were provided for scalar gaps
and attractive interaction kernels in
\cite{billard1968existence,yang1991bardeen}. In the more recent
literature, the scalar BCS gap equation within the continuum
approximation, neglecting the lattice structure, was analyzed in
\cite{hainzl2008critical}, obtaining bounds on the critical
temperature as well as on the gap, and results on the critical
scaling. In \cite{hainzl2008bcs}, the continuum BCS free energy
functional for general pair interactions is analysed, which was
derived by Leggett in \cite{leggett1980diatomic} based on the
original work by Bardeen \cite{bardeen1957theory}. Symmetries of
minimisers are analysed in \cite{deuchert2018persistence}. The
extension to multiple bands is discussed in \cite{yang2005mathematical}.
For a review on mathematical results concerning the BCS gap equation,
see \cite{hainzl2016bardeen}. Only few works exist that treat
numerical analysis aspects of the BCS equation and of related models.
A convergence analysis comparing self-consistent iterations with a
direct minimisation approach for models of Hartree-Fock type has been
provided in \cite{cances2021convergence}; its setting of
density matrices as projectors of fixed rank does not include pairing.
Here, we note that self-consistent iteration is the
standard approach for solving the gap equation in the physics
literature for attractive interaction kernels. Numerical analysis
aspects of Hartree-Fock-Bogoliubov theory, a generalisation of the
BCS model, are discussed within a continuum approximation in
\cite{lewin2014numerical}. The numerical solution for large scale
Hartree-Fock-Bogoliubov systems for finite lattice systems has been
recently treated in \cite{lin2021numerical}.

In this work, we provide an efficient and precise method for
computing the superconducting gap function for unconventional
superconductors at zero or finite temperature within the BCS
framework. The method applies to any Bravais lattice, any
non-repulsive interaction composed of power-laws plus short-ranged
contributions, any dispersion relation, and both singlet and triplet
channels. In particular, it includes interaction kernels that are not
absolutely summable, important cases being the 3D Coulomb interaction
in a one-dimensional lattice and the isotropic dipole interaction for
three-dimensional lattices. We first demonstrate that solutions of
the gap equation correspond to stationary points of a free energy
functional. We show that this functional is always bounded from below
and that self-consistent iteration, an often used method in the
physics literature, never increases the free energy and leads to a
vanishing residual of the gap equation. After rewriting the gap
equation in a two-field formulation in terms of the gap and the
anomalous density, we then show that a Galerkin B-spline
discretisation in momentum space allows for an explicit
representation of the singular convolutions in terms of rapidly
convergent lattice sums. It leads to a block-circulant Galerkin
matrix, whose eigenvalues can be efficiently evaluated. Thus, any
iteration is log-linear in numerical work and linear in memory. For
the integer lattice at zero temperature and power-law interactions,
we subsequently derive new analytical singlet and triplet solutions
of the gap equation, based on engineering a suitable dispersion for a
prescribed gap. These solutions have reduced regularity and allow for
a reliable benchmarking of our algorithm. After testing the method
against the analytical solutions as well as against highly resolved
numerical data, we use the method to study nodal Weyl superconductors
on the three-dimensional face-centered-cubic lattice.

This work is structured as follows. In Sec.~\ref{sec:model}, we
introduce some necessary notation, define the gap equation and
introduce the BCS energy functional, and present the theoretical results. In
Sec.~\ref{sec:galerkin}, we then formulate the Galerkin B-spline
discretisation and derive the eigenvalues of the arising matrices.
The analytical solutions to the gap equation are provided in
Sec.~\ref{sec:analytical}. Detailed numerical benchmarks and a
convergence analysis are provided in Sec.~\ref{sec:benchmarks}.
Finally, we display solutions for nodal 3D superconductors with
long-range interactions in Sec.~\ref{sec:experiments}. We draw our
conclusions and provide an outlook in Sec.~\ref{sec:outlook}

\section{Gap Equation, Energy Functional, and Main Properties}
\label{sec:model}


We will begin our discussion with the definition of mono-atomic
Bravais lattices. These were first introduced by the French
mathematician Auguste Bravais in 1850, \cite{BRAVAIS:1850}.
\begin{definition}[Bravais lattices]
  Let $A\in \mathds R^{d\times d}$ invertible. We call the periodic
  arrangement of points $\Lambda=A\mathds Z^d$ a mono-atomic lattice
  or Bravais lattice. Its elementary cell $E$ is defined as
  $E=A\mathds T^d$ with the torus $\mathds T=[0,1)$. The reciprocal
  lattice is defined as $\Lambda^\ast = A^{-\top} \mathds Z^d$ with
  reciprocal elementary cell $E^\ast = A^{-\top}\mathds T^d$.
\end{definition}

In the following, we will make use of the reciprocal lattice cell
$E^\ast$ as the integration domain in integrals over
$\Lambda^\ast$-periodic functions. This cell is closely related to
the Brillouin zone, used in the physics literature, which is defined as
\[
  \mathrm{BZ} = \big\{\bm p\in \mathds R^d: |\bm p|\le |\bm p-\bm
  q|,\quad \bm q\in \Lambda^\ast\setminus\{\bm 0\}\big\}.
\]
In particular, $|\mathrm{BZ}|=|E^\ast|$ and integrals of
$\Lambda^\ast$-periodic functions over $\mathrm{BZ}$ and $E^\ast$
coincide, as the Brillouin zone can be seen as a centered
rearrangement of parts of the reciprocal lattice cell, which exhibits
the same symmetries as the underlying lattice. In this work, we use
the reciprocal lattice cell as the computational cell and display
numerical results within the Brillouin zone.

In the following, we require  Fourier series on lattices.

\begin{definition}[Lattice Fourier series]
  For a Bravais lattice $\Lambda=A\mathds Z^d$, the lattice Fourier
  transform of $f\in \ell^1(\Lambda)$ is defined as
  \[
    \hat f(\bm p) = \mathcal F_\Lambda [f](\bm p)= \sum_{\bm x\in
    \Lambda} e^{-\imath 2\pi\, \bm x\cdot \bm p} f(\bm x)\,,\quad \bm
    p\in E^\ast,
  \]
  which we extend to a $\Lambda^\ast$-periodic function on $\mathds
  R^d$. The inverse transform, for $g\in L^1(E^\ast)$, is defined as
  \[
    \mathcal F_\Lambda^{-1}[g](\bm x) =
    \frac{1}{|E^\ast|}\int\limits_{E^\ast} e^{\imath 2\pi\,\bm x\cdot
    \bm p}g(\bm p)\,\mathrm d \bm p\,,\quad \bm x\in \Lambda.
  \]
  Finally, the periodic convolution on $E^\ast$ is defined for $\hat
  f,\hat g\in L^2(E^\ast)$ as
  \[
    (\hat f\ast \hat g)(\bm p)=\frac{1}{|E^\ast|}\int\limits_{E^\ast}
    \hat f(\bm p-\bm q)\hat g(\bm q)\, \mathrm{d}\bm q\,.
  \]
\end{definition}

We consider superconductivity with Cooper pairs of zero total
momentum, i.e.\ electrons pairing at wave vectors\footnote{The wave
  vector is often referred to as the momentum in the physics
literature, as they are related by the de Broglie relation.}
$\boldsymbol{p}$ and $-\boldsymbol{p}$ with $\bm p\in E^\ast$, for a
tight-binding model on a Bravais lattice $\Lambda\subseteq \mathds
R^d$ and a spin-independent electron-electron interaction kernel
$K\in \ell^1(\Lambda,\mathds R)$ within the mean-field approximation,
focusing in the numerical treatment on attractive (in our sign
convention positive) kernels. Tunneling of electrons between lattice
sites is described by the spin-independent and even hopping amplitude
$\tau:\Lambda\to \mathds R$ (which includes the Planck constant
$\hbar$ and is in units of energy), which is often compactly
supported. The dispersion relation $\xi(\bm p)$ is defined as the
lattice Fourier transform of the hopping amplitude shifted by the
chemical potential $\mu\in \mathds R$,
\[
  \xi(\bm p)= \mathcal F_\Lambda[\tau](\bm p)-\mu\,.
\]
Then $\xi$ is real due to evenness, and analytic on $E^\ast$ if
$\tau$ is compactly supported. If $\tau$ extends over the whole
lattice, $\xi$ loses regularity. In the following, we will assume
only that $\xi$ is even and continuous, corresponding to $\tau \in
\ell^1(\Lambda,\mathds R)$ even. For a review of microscopic
unconventional superconductivity in the tight-binding limit, see
\cite{sigrist1991phenomenological}.

A superconducting gap function $\hat \Delta$ in the momentum space
lies in the following function space. The additional symmetry
constraint results from fermionic anticommutation rules.
\begin{definition}[Superconducting gap function, solution space, and
  inner product]
  \label{def:gap_space}
  For a Bravais lattice $\Lambda=A\mathds Z^d$ with reciprocal
  lattice cell $E^\ast$, let
  \[
    X=\big\{ \hat \Delta\in L^2(E^\ast,\mathds C^{2\times 2}):  \hat
    \Delta^\top(\bm p)=-\hat \Delta(-\bm p),~\bm p\in E^\ast\big\}\,,
  \]
  We call $\hat \Delta\in X$ the momentum-space superconducting gap
  matrix. Equivalently, the real space gap lies in $\Delta\in
  \mathcal F_\Lambda^{-1} X$, thus $\Delta \in \ell^2(\Lambda,\mathds
  C^{2\times 2})$ obeying $\Delta^\top(\bm x)=-\Delta(-\bm x)$ for
  all $\bm x\in \Lambda$. For general matrix valued functions $f,g\in
  L^2(E^\ast,\mathds C^{2\times 2})$ (with no restrictions on
  symmetry), we further define the inner product
  \[
    \langle f,g \rangle = \frac{1}{|E^\ast|} \int\limits_{E^\ast}
    \mathrm{Tr} \big(f(\bm p)^\ast g(\bm p) \big)\,\mathrm d \bm p\,.
  \]
\end{definition}

In order to write down the BCS gap equation, we need to introduce the
Fermi surface as the zero-set of $\xi$, as well as two operators
acting on the gap function $\hat \Delta$ in momentum space, namely
the quasi-particle energy operator $E[\hat \Delta]$ and the anomalous
density operator $G[\hat \Delta]$.
\begin{definition}[Fermi-surface, quasi-particle energy, and anomalous density]
  \label{def:anomalous_density}
  Consider $\xi\in C(E^\ast,\mathds R)$ even with zero-set
  \[
    \mathcal C=\{\bm p\in E^\ast: \xi(\bm p)=0\}\,.
  \]
  We call $\mathcal C$ the Fermi-surface, which we assume to have
  zero measure. For $\hat \Delta\in X$, we define the self-adjoint
  quasi-particle energy operator $E[\hat \Delta]$ as
  \[
    E[\hat \Delta]=\big(\xi^2\mathds 1_2+\hat \Delta^\ast\hat \Delta\big)^{1/2},
  \]
  where the right-hand side is to be understood in terms of the matrix power.
  Further, for coldness $\beta=1/(k_B T)\in (0,\infty]$
  for $k_B$ the Boltzmann
  constant and $T\ge 0$ the temperature, we define the anomalous
  density $G_\beta:X\to X$,
  \[
    G_\beta[\hat \Delta]=\frac{1}{2}\hat \Delta \big(E[\hat
    \Delta]\big)^{-1} \tanh\big(\beta E[\hat \Delta]/2\big).
  \]
  Note that $\tanh\big(\beta E[\hat \Delta]/2\big)=(\mathds
    1_2-2N_\beta[\hat \Delta]
  )$
  with $N_\beta[\hat \Delta]$ the quasi-particle occupancy matrix given as
  \[
    N_\beta[\hat \Delta]= (1+e^{\beta E[\hat \Delta]})^{-1}.
  \]
  The case $T=0$, corresponding to $\beta\to \infty$, is obtained in
  the sense of the limit, where the $\tanh$ term is replaced by unity.
\end{definition}

The well-definedness of the anomalous density operator and its
properties, in particular the fact that it preserves the symmetry of
$\hat \Delta$, are discussed in Lem.~\ref{lem:aux_lemma_gb} in the appendix.

We can now state the Bardeen--Cooper--Schrieffer gap equation as a
nonlinear matrix-valued convolution equation of Hammerstein type.

\begin{definition}[Bardeen-Cooper-Schrieffer gap equation]
  \label{def:gap_eq}
  Let $\hat \Delta\in X$ and consider an even interaction kernel
  $\hat K\in L^2(E^\ast,\mathds R$). The Bardeen--Cooper--Schrieffer
  gap equation reads
  \begin{equation*}
    \hat \Delta=\hat K\ast  G_\beta[\hat \Delta]
  \end{equation*}
  where the convolution operator is applied component-wise.
\end{definition}

The gap equation is the Euler-Lagrange equation of the free energy
$F_\beta[\hat \Delta]$, see
\cite{leggett1980diatomic,fanfarillo2025superconductivity,hainzl2008bcs}.
We formulate the free energy for triplet superconductors on a Bravais
lattice $\Lambda$ as follows.

\begin{definition}[BCS free energy functional]
  \label{def:functional}
  Let $\hat \Delta \in X$, $\hat K\in L^2(E^\ast,\mathds R)$, and
  $\beta\in (0,\infty]$. The BCS free energy functional $F_\beta:X\to
  \mathds R$ reads
  \begin{equation*}
    F_\beta[\hat \Delta]=E_{\mathrm{pot},\beta}[\hat
    \Delta]+E_{\mathrm{kin},\beta}[\hat \Delta]-T S[\hat \Delta]\,,
  \end{equation*}
  with $E_{\mathrm{pot},\beta}$ the potential energy,
  $E_{\mathrm{kin},\beta}$ the kinetic energy, and $S$ the entropy,
  all measured per lattice site. Defining the momentum-averaged trace as
  \[
    \mathrm{tr}\,F =\langle \mathds 1_2, F\rangle =
    \frac{1}{|E^\ast|} \int\limits_{E^\ast} \mathrm{Tr}\, F(\bm
    p)\,\mathrm d \bm p\,,
  \]
  the three contributions read
  \begin{align*}
    E_{\mathrm{pot},\beta}[\hat \Delta]&=-\frac{1}{2}\langle
    G_\beta,\hat K\ast G_\beta\rangle\,,\\
    E_{\mathrm{kin},\beta}[\hat \Delta]&=
    \frac{1}{2}\mathrm{tr}\bigg(\xi \mathds 1_2-\xi^2 E[\hat
    \Delta]^{-1}(\mathds 1_2-2 N_\beta[\hat \Delta])\bigg)\,, \\
    -TS[\hat \Delta]&=\frac{1}{\beta} \mathrm{tr}\Big(N_\beta[\hat
      \Delta] \log N_\beta[\hat \Delta]+(1-N_\beta[\hat
    \Delta])\log(1-N_\beta[\hat \Delta])\Big)\,,
  \end{align*}
  The zero-temperature functional is defined in terms of the limit
  $\beta\to \infty$.
\end{definition}

In the following lemma, we combine kinetic energy and entropy into a
single potential. Importantly, the resulting functional only depends
on $\hat \Delta^\ast \hat\Delta$, or, equivalently, on $E[\hat \Delta]$,
which will be relevant for the subsequent analysis.
\begin{lemma}
  \label{lem:energy_reformulation}
  Let $\hat \Delta\in X$ and $0<\beta<\infty$. Then
  \begin{equation*}
    F_{\mathrm{kin},\beta}[\hat \Delta]=E_{\mathrm{kin},\beta}[\hat
    \Delta]-TS[\hat\Delta] = \frac{1}{2}\mathrm{tr}\big(\xi\mathds
    1_2-\phi_\beta(E[\hat \Delta])\big)+\langle \hat \Delta,
    G_\beta[\hat \Delta] \rangle,
  \end{equation*}
  with the real-analytic and even function $\phi_\beta:\mathds R\to \mathds R$,
  \[
    \phi_\beta(t)=\frac{2}{\beta}\log\Big(2\cosh(\beta t/2)\Big).
  \]
  In the limit $\beta\to\infty$, we obtain $\phi_\infty(t)=|t|$,
  hence $\phi_\infty(E[\hat \Delta])=E[\hat \Delta]$.
\end{lemma}
\begin{proof}
  We first rewrite the entropy contribution. With
  \begin{align*}
    N_\beta\log N_\beta &=-N_\beta\beta E -N_\beta\log(\mathds
    1_2+e^{-\beta E}),
  \end{align*}
  we have
  \[
    -TS=\mathrm{tr}\big(-E N_\beta-\tfrac{1}{\beta}\log(\mathds
    1_2+e^{-\beta E})\big)=\mathrm{tr}\big(-E N_\beta
    +\tfrac{1}{2}E-\tfrac{1}{2}\phi_\beta(E)\big).
  \]
  In the kinetic part, expand
  \[
    \xi^2 E^{-1}=(E^2-\hat \Delta^\ast \hat \Delta) E^{-1} =E-\hat
    \Delta^\ast \hat \Delta E^{-1},
  \]
  yielding
  \[
    E_\mathrm{kin,\beta}= \mathrm{tr}\big(\tfrac{1}{2}\xi\mathds
    1_2+E N_\beta -\tfrac{1}{2}E\big)+\langle \hat
    \Delta,G_\beta[\hat \Delta]\rangle.
  \]
  When combining the two contributions, the terms $\pm E N_\beta$ and
  $\pm \tfrac{1}{2}E$ cancel, and the statement follows.
\end{proof}

The functional from Def.~\ref{def:functional} arises when applying a
particular trial wave function, called a quasi-free translational
invariant state, to the many-body Hamiltonian and computing the
energy per particle and subsequently adding the entropy contribution
of the state. The functional is therefore an upper bound to an
expectation value of a self-adjoint operator that is bounded from
below. For this reason, the functional has to be bounded from below
as well. We offer a rigorous proof in the following proposition. For
a treatment of the functional within a continuum approximation,
discarding the microscopic lattice structure, see \cite{hainzl2008bcs}.
\begin{proposition}[Boundedness of the BCS functional]
  \label{prop:boundedness}
  For $\hat \Delta\in X$, the functional from Definition
  \ref{def:functional} is bounded from below. Setting
  \[
    K_\mathrm{sup}=\sup_{\bm x\in \Lambda} K(\bm x)\,,
  \]
  the bound reads
  \[
    F_\beta[\hat \Delta]\ge F_\beta[0]-\tfrac{1}{4}K_\mathrm{sup}\,.
  \]
\end{proposition}
\begin{proof}
  As $K$ is bounded, $K_\mathrm{sup}$ exists. Separate
  $F_\beta=E_{\mathrm{pot,\beta}}+ F_{\mathrm{kin},\beta}$. For the
  bound on the potential energy, we first write the potential in real space form
  \[
    E_\mathrm{pot,\beta}[\hat \Delta]=-\frac{1}{2}\sum_{\bm x\in
    \Lambda} K(\bm x) \Vert g_\beta[\hat \Delta](\bm x)\Vert^2_F
  \]
  with $g_\beta[\hat \Delta]=\mathcal F_\Lambda^{-1} G_\beta[\hat
  \Delta]$ and $\Vert \cdot \Vert_F$ the Frobenius norm.
  Subsequently, bound $K$ by its supremum and apply Parseval's theorem to obtain
  \[
    E_\mathrm{pot,\beta}[\hat \Delta] \ge -\frac{K_\mathrm{sup}}{2}
    \langle G_\beta[\hat \Delta], G_\beta[\hat \Delta]\rangle.
  \]
  For a singular value $\sigma$ of $\hat \Delta$ at given wave
  vector, the eigenvalue of $G_\beta^\ast[\hat \Delta] G_\beta[\hat
  \Delta]$ reads
  \[
    \frac{\sigma^2}{4 e^2}\tanh{(\beta e/2)}^2\le \frac{1}{4},\quad
    e=\sqrt{\xi^2+\sigma^2}.
  \]
  For the sum of the two eigenvalues and the momentum average, we
  therefore obtain
  \[
    \langle G_\beta[\hat \Delta], G_\beta[\hat \Delta]\rangle\le
    \frac{1}{2},\quad \mathrm{thus}\quad E_\mathrm{pot,\beta}[\hat
    \Delta]\ge -\frac{K_\mathrm{sup}}{4}.
  \]
  Note that $K_\mathrm{sup}\ge 0$, as Riemann--Lebesgue implies that
  $K(\bm x)$ vanishes as $\vert \bm x\vert\to \infty$.

  For the bound on $F_{\mathrm{kin},\beta}$, notice that
  $F_{\mathrm{kin},\beta}[\hat \Delta]=\mathrm{tr}\,f_\xi(E[\hat
  \Delta])$, where at each wave vector the scalar function
  $f_\xi:[\,|\xi|,\infty)\to\mathds R$,
  \[
    f_\xi(e)=\frac12\big(\xi-\phi_\beta(e)\big)+\frac{e^2-\xi^2}{2e}\phi_\beta'(e),
  \]
  is applied to the eigenvalues of $E[\hat \Delta]$.
  For a singular value $\sigma$ of $\hat \Delta$, the eigenvalues of
  $E[\hat \Delta]$ obey $e=\sqrt{\xi^2+\sigma^2}\ge |\xi|$ with
  equality only for $\sigma=0$. It is therefore sufficient to show
  that $f_\xi(e)$ is strictly increasing for $e> |\xi|$ to show that
  the minimum of $F_{\mathrm{kin},\beta}$ is obtained at $\hat
  \Delta=0$. The derivative of $f_\xi$ reads
  \[
    f_\xi'(e)=\frac{\xi^2}{2e^2}\,\phi_\beta'(e)+\frac{e^2-\xi^2}{2e}\,\phi_\beta''(e)\ge
    0,\quad e> |\xi|,
  \]
  as $\phi_\beta$ is both increasing on the non-negative reals and convex with
  \[
    \phi_\beta'(e)=\tanh(\beta e/2)\ge 0,\quad
    \phi_\beta''(e)=\tfrac{\beta}{2}\big(1-\tanh^2(\beta e/2)\big)>0\,,
  \]
  Thus, the minimum $f_\xi$ is attained at $e=|\xi|$. This yields
  that $F_{\mathrm{kin},\beta}[\hat \Delta]$ is minimal at $\hat
  \Delta =0$. Furthermore, as the potential energy term vanishes for
  $\hat \Delta =0$, we have that
  \[
    F_{\mathrm{kin},\beta}[\hat \Delta]\ge \mathrm{tr}
    f_\xi\big(|\xi|\mathds 1_2\big)= F_{\mathrm{kin},\beta}[0]=F_\beta[0]\,.
  \]
  Combining the bounds yields the statement.
\end{proof}

The above bound shows that superconductivity can only arise, when $K$
is positive at some displacement vector $\bm x\in \Lambda$,
meaning that
electrons separated by this displacement vector attract each other. In
the absence of such attraction, the normal state $\hat \Delta=0$
forms the global energy minimum.

Finally, we show that stationary points of the BCS functional are
solutions to the gap equation.
\begin{theorem}
  \label{thm:gap_eq_is_stationary_point}
  Let $0<\beta<\infty$ and $\hat \Delta\in X$.  Then $\hat \Delta$ is
  a stationary point of $F_\beta$ if and only if $\hat \Delta$ solves
  the BCS gap equation from Def.~\ref{def:gap_eq}.
\end{theorem}
\begin{proof}
  Let $\hat H\in X$ be an arbitrary perturbation direction and let
  $\delta$ denote the derivative at $\varepsilon=0$ along $\hat
  \Delta +\varepsilon \hat H$,
  thus
  \[
    \delta F_\beta = \frac{\partial}{\partial \varepsilon}
    F_\beta[\hat \Delta+\varepsilon\hat H]\Big\vert_{\varepsilon=0}\,,
  \]
  and let $\delta G_\beta[\hat \Delta]$ denote the pointwise
  derivative of $\varepsilon\mapsto G_\beta[\hat \Delta+\varepsilon\hat H]$.
  The residual $R$ of the BCS gap equation reads
  \[
    R=\hat \Delta-\hat K\ast G_\beta[\hat \Delta]\,.
  \]
  We first show that the gap equation is well-defined in $X$. As
  $G_\beta[\hat \Delta]\in X$ by Lem.~\ref{lem:aux_lemma_gb} and
  convolution with the even kernel $\hat K$ preserves the symmetry,
  we have $R\in X$.

  By the reformulation of Lem.~\ref{lem:energy_reformulation}, the
  free energy reads
  \[
    F_\beta[\hat \Delta]= \tfrac{1}{2}\big\langle \mathds
    1_2,\xi\mathds 1_2-\phi_\beta(E[\hat \Delta])\big\rangle+\langle
    \hat \Delta, G_\beta[\hat \Delta] \rangle-\tfrac{1}{2}\langle
    G_\beta[\hat \Delta],\hat K\ast G_\beta[\hat \Delta]\rangle,
  \]
  where we recall that
  \[
    \phi_\beta(t)=\frac{2}{\beta}\log\Big(2\cosh(\beta t/2)\Big).
  \]
  We begin with the variation of the first term in the free energy.
  Since $\phi_\beta$ is even and real analytic,
  $\varphi_\beta(t)=\phi_\beta(\sqrt t)$ is smooth on $[0,\infty)$ with
  \[
    \varphi_\beta'(t) = \frac{1}{2}\tanh(\beta \sqrt{t}/2)\frac{1}{\sqrt{t}},
  \]
  where we replace $t$ by $E[\hat \Delta+\varepsilon \hat H]^2$ in
  the following. Then notice that
  \[
    E[\hat \Delta+\varepsilon \hat H]^2=\xi^2\mathds 1_2 +(\hat
    \Delta^\ast+\varepsilon \hat H^\ast)(\hat \Delta+\varepsilon \hat H)
  \]
  is a quadratic polynomial in $\varepsilon$ with variation
  \[
    \delta \big( E[\hat \Delta]^2\big)=\hat \Delta^\ast \hat H+\hat
    H^\ast \hat \Delta .
  \]
  By the formula for the derivative of a matrix function under a
  trace, see \cite[Thm.~3.11]{Higham2008}, we then find
  \[
    \delta (\tfrac{1}{2}\langle \mathds 1_2, -\phi_\beta(E[\hat
    \Delta]) \rangle) =-\tfrac{1}{2}\big\langle \mathds
    1_2,\tfrac{1}{2}E[\hat \Delta]^{-1}\tanh(\beta E[\hat
    \Delta]/2)\big(\hat \Delta^\ast \hat H+ \hat H^\ast\hat
    \Delta\big)\big\rangle.
  \]
  After cyclic reordering within the trace, inserting the definition
  of $G_\beta[\hat \Delta]$ yields
  \[
    \delta (\tfrac{1}{2}\langle \mathds 1_2, -\phi_\beta(E[\hat
    \Delta]) \rangle) =-\mathrm{Re}\langle \hat H,G_\beta[\hat \Delta]\rangle.
  \]

  For the second term in the energy, we use that $\Delta^\ast
  G_\beta[\hat \Delta]$ is self-adjoint and thus $\langle \hat
  \Delta, G_\beta[\hat \Delta] \rangle $ is real. This yields
  \[
    \delta (\langle\hat \Delta, G_\beta[\hat \Delta] \rangle)=
    \mathrm{Re}\langle \hat H,G_\beta[\hat \Delta]\rangle +
    \mathrm{Re}\langle \hat \Delta,\delta G_\beta[\hat \Delta] \rangle.
  \]
  Hence, the variations of the first and second term in the energy combine into
  \[
    \mathrm{Re}\langle \hat \Delta,\delta G_\beta[\hat \Delta] \rangle.
  \]
  Finally, for the remaining potential energy part, we have
  \[
    -\frac{1}{2}\delta \langle G_\beta[\hat \Delta],\hat K\ast
    G_\beta[\hat \Delta]\rangle=-\frac{1}{2}\Big(\langle \delta
      G_\beta[\hat \Delta],\hat K\ast G_\beta[\hat
      \Delta]\rangle+\langle G_\beta[\hat \Delta],\hat K\ast \delta
    G_\beta[\hat \Delta]\rangle\Big).
  \]
  Now as the convolution with the even function $\hat K$ is
  self-adjoint on $L^2(E^\ast,\mathds C^{2\times 2})$, the right hand
  side can be combined into
  \[
    \mathrm{Re}\langle - \hat K\ast G_\beta[\hat \Delta],\delta
    G_\beta[\hat \Delta] \rangle.
  \]
  Note that all derivatives exist by dominated convergence. Given
  $|\varepsilon|\le 1$, Lem.~\ref{lem:aux_lemma_gb}
  bounds the $\varepsilon$-derivatives of the integrands by \[
    c\,(1+\Vert\hat \Delta\Vert_F+\Vert\hat H\Vert_F)\Vert\hat
  H\Vert_F\]  for some $c>0$, and the difference quotients of
  $G_\beta[\hat \Delta+\varepsilon\hat H]$ by
  $\tfrac{\beta}{4}\Vert\hat H\Vert_F$, pointwise in $\bm p$.

  Combining all three contributions yields
  \[
    \delta F_\beta[\hat \Delta]=\mathrm{Re}\,\langle \hat \Delta-\hat
    K\ast G_\beta[\hat \Delta],\delta G_\beta[\hat \Delta]
    \rangle=\mathrm{Re}\,\langle R,\delta G_\beta[\hat \Delta]\rangle\,,
  \]
  with $R$ the residual of the gap equation. Hence, all solutions of
  the gap equation are stationary points of the functional.

  For the reverse direction, let $\hat \Delta$ be a stationary point
  of $F_\beta$, thus $\delta F_\beta[\hat \Delta]$ is zero along
  every direction $\hat H\in X$. Choose the gap equation residual as
  the direction, $\hat H=R$, which is possible as $R\in X$. Then
  \[
    \delta F_\beta[\hat
    \Delta]=\frac{1}{|E^\ast|}\int\limits_{E^\ast}\mathrm{Re}\,
    \mathrm{Tr}\big(R(\bm p)^\ast \delta G_\beta[\hat \Delta](\bm
    p)\big)\,\mathrm d \bm p = 0\,.
  \]
  Now, by Lem.~\ref{lem:aux_lemma_gb}.2, the integrand is positive
  whenever $R(\bm p)\neq 0$ and is zero otherwise. Thus, we obtain
  that $R=0$ almost everywhere, which is the gap equation.
\end{proof}

It is important to realise that only solutions of the gap equation
that correspond to local minima of the free energy functional are
stable under infinitesimal perturbation and, therefore, physically
meaningful. It is common practice in the physical literature to solve
the gap equation through fixed-point iteration, where the original
energy functional has been discarded or, in many cases, forgotten.
The following Proposition shows that fixed-point iteration of the gap
equation never increases the energy in each iteration step and drives
the residual of the gap equation to zero for non-negative kernels
$K$. Meanwhile, one can show that stationary points that are not
local minima are unstable fixed points of the iteration. Up to a set
of measure zero of starting points, iteration thus converges to local
minima. In other words, solving the gap equation under fixed-point
iteration is equivalent to preconditioned gradient descent for the
free energy, explaining why the common practice is meaningful for
non-repulsive interaction kernels. It is an instance of the
concave-convex procedure from \cite{yuille2003concave}. For a
comparison of self-consistent iteration and direct minimisation for
Hartree-Fock and Kohn-Sham models, see \cite{cances2021convergence}.
In our setting, the fixed
dispersion and sign of $K$ guarantee that the residual converges to
zero. Note that in the Hartree-Fock-Bogoliubov case, the iteration
may oscillate between states \cite{lewin2014numerical}.

\begin{proposition}[Convergence under iteration]
  \label{prop:convergence}
  Let $0<\beta<\infty$ and $K\ge 0$ even. Choose an initial gap $\hat
  \Delta_0\in X$ and consider the iteration
  \[
    \hat \Delta_{j+1} = \hat K\ast G_j,\quad G_j=G_\beta[\hat \Delta_j],
  \]
  with the residual $R_j=\hat \Delta_j-\hat K\ast G_j$ at the $j$-th iteration.
  It then holds that
  \[
    F_\beta[\hat \Delta_{j+1}]-F_\beta[\hat \Delta_j]\le
    -\frac{1}{2}\langle G_{j+1}-G_j,\hat K\ast (G_{j+1}-G_j) \rangle\le 0,
  \]
  and
  \[
    \sum_{j=1}^\infty \Vert R_j\Vert^2\le 2 K_\mathrm{sup}
    \Big(F_\beta[\hat \Delta_0]-F_\beta[0]+\tfrac{1}{4}K_\mathrm{sup}\Big).
  \]
  The free energy therefore decreases monotonically and the residual
  of the gap equation converges to zero as $j\to\infty$.
\end{proposition}

The proof relies on properties of the anomalous density, seen as a
pointwise function of $\hat \Delta(\bm p)$, from
Lem.~\ref{lem:aux_lemma_gb} and is given in the appendix.

The definitions of the gap equation and of the free energy
functional, as well as Prop.~\ref{prop:boundedness},
Thm.~\ref{thm:gap_eq_is_stationary_point}, and
Prop.~\ref{prop:convergence}, extend without modification to any
kernel $\hat K\in L^1(E^\ast)$, with the lattice Fourier transform
understood in the sense of tempered distributions. In particular, the
kernel only enters through its Fourier coefficients $K(\bm x)$, which
are bounded and vanish at infinity. In particular, this includes all
regularised power laws $\mathcal K_\nu:\Lambda\to \mathds R$, defined
as  $\mathcal K_\nu(\bm x)=\vert \bm x\vert^{-\nu}$, for $\bm x\neq
\bm 0$, and $\mathcal K_\nu(\bm 0)=0$, under the restriction $\nu>0$.
Here, absolute summability is lost, yet integrability of the lattice
Fourier transform is sufficient for all properties discussed above to hold.
The lattice Fourier transform of these kernels yields the Epstein
zeta function \cite{epstein1903theorieI,epstein1903theorieII}, which
we define below.

\begin{definition}[Epstein zeta function]
  \label{def:epstein}
  For a Bravais lattice $\Lambda$, and $\bm x,\bm p \in \mathds R^d$
  and $\nu\in \mathds C$, the Epstein zeta function is defined as
  \[
    Z_{\Lambda,\nu}(\bm x,\bm p) = \sum_{\bm z\in \Lambda} e^{-\imath
    2\pi\,\bm z\cdot \bm p} \mathcal K_\nu(\bm z-\bm x)=\mathcal
    F_\Lambda\big[\mathcal K_\nu(\cdot -\bm x)\big](\bm p),\quad
    \mathrm{Re}(\nu)>d,
  \]
  which can be meromorphically continued to $\nu\in \mathds C$. In
  case that $\bm x=\bm 0$, we adopt the shortened notation
  $Z_{\Lambda,\nu}(\bm p)$.
\end{definition}

Its analytical properties and the efficient numerical evaluation of
the Epstein zeta function were recently described in
\cite{buchheit2024epstein}, together with an implementation in the
high-performance C-library EpsteinLib, which is publicly available
\cite{epsteinlib}. A recent extension to iterated lattice sums with a
graph structure was provided in
\cite{buchheit2025epstein_method,buchheit2026graph,gzl2026}, yielding
access to efficient series expansions of gapped quantum lattice
models \cite{duft2026}.

In the following, we will focus on positive\footnote{In our sign
  convention, this corresponds to purely attractive interactions
between electrons.} interaction kernels of the form
\begin{equation}
  \label{eq:interaction}
  K(\bm x) = a(\bm x) + \sum_{j=1}^n b_j \mathcal K_{\nu_j}(\bm x)\,,
\end{equation}
with $a:\Lambda\to \mathds R$ positive, even, and compactly
supported, $b_j>0$, and $\nu_j>0$. Their lattice Fourier transform
can then be written in terms of Epstein zeta functions,
\[
  \hat K(\bm p) = \hat a(\bm p) +\sum_{j=1}^n b_j Z_{\Lambda,\nu_j}(\bm p)\,.
\]
Here, the compactly supported part of the interaction describes the
near-field, whereas long-range tails are efficiently describable
through different exponents $\nu_j$. The inclusion of anisotropic
interactions is possible via the recent anisotropic generalisation of
the Epstein zeta function in \cite{buchheit2026aniso}.

From now on, we will adopt dimensionless units, writing positions in
units of a reference length $a_0$, for instance the smallest distance
between lattice points, wave vectors in units of $1/a_0$, and
energies in units of a reference band energy $e_0$, e.g.~the maximum
absolute value of the hopping amplitude.

\section{Galerkin Approximation}
\label{sec:galerkin}
In this section, we provide a Galerkin discretisation in terms of
periodic B-splines of the gap equation on the reciprocal lattice cell
$E^\ast$ from Def.~\ref{def:gap_eq}. The method is applicable for any
coldness $\beta\in (0,\infty]$, any Bravais lattice $\Lambda=A\mathds
Z^d$, with $A\in \mathds R^{d\times d}$ invertible, any even
interaction kernel $\hat K\in L^1(E^\ast,\mathds R)$, in particular
sums of power-laws in the form of Eq.~\eqref{eq:interaction} with
exponents $\nu_j>0$ plus short-ranged parts, and both singlet and
triplet channels. The method was first outlined in
\cite{buchheit2026numerical} for square lattices at $\beta=\infty$
and single power-law interactions. For the iteration, we currently
restrict the method to $K\ge 0$, while the remaining parts apply to
general kernels.

The method rests on two ingredients. First, we apply a two-field
formulation, splitting the gap equation into a linear convolution
equation and a pointwise nonlinearity. The second observation is that
the matrix of the convolution operator after B-spline discretisation
is block-circulant with eigenvalues given by explicit lattice sums
involving $K(\bm x)$, $\bm x\in \Lambda$. The section culminates in
an explicit expression for these eigenvalues, reducing, once
precomputed, every application of the convolution operator to fast
Fourier transforms. By allowing efficient and precise evaluations of
the singular convolutions at log-linear cost, it allows for the
simulation of three-dimensional superconductors with complicated
singularity structures within minutes on standard desktop hardware,
as demonstrated in Sec.~\ref{sec:experiments}.

\subsection{Two-Field Formulation}
We begin by writing the solution of the gap equation in terms of two
unknown fields, the gap function $\hat \Delta\in X$ and the anomalous
density $G\in X$, with
\begin{equation}
  \label{eq:two-field-formulation}
  \hat \Delta =\hat K\ast G,\quad G=G_\beta[\hat \Delta].
\end{equation}
The first equation is linear, namely a periodic convolution with the
kernel $\hat K$, while the second is nonlinear, yet only depends
pointwise on the values of $\hat \Delta(\bm p)$ and $\xi(\bm p)$. In
our numerical B-spline discretisation of the gap equation, we use
this formulation
with two fields as it separates the linear periodic convolution in
the first equation from the pointwise nonlinear relation between $G$
and $ \hat \Delta$ in the second equation.

Note that the normal state $\hat \Delta=0$, $G=0$, is always a
solution. In general, if $\hat \Delta$ solves the gap equation, then
for any constant unitary matrix $V\in \mathds C^{2\times 2}$ also the
pair $V\hat \Delta V^\top$, $V G V^\top$ solves
Eq.~\eqref{eq:two-field-formulation}. For $V=e^{\imath
\varphi}\mathds 1_2$, $\varphi \in \mathds R$, this corresponds to a
global phase, while for singlet solutions
\[
  \hat \Delta=
  \begin{pmatrix}
    0 & \psi \\
    -\psi & 0
  \end{pmatrix},
\]
the family also reduces to a phase factor for $\psi$. An iterative
solver then converges to one representative of such a family, each of
them being physically equivalent.

\subsection{Fourier Series, Sobolev Spaces, and Convolution Operators}
In this section, we first introduce some basic notation, briefly
discuss Fourier series and their connection to Sobolev spaces, and
then discuss implications for the convolution operator from
Eq.~\eqref{eq:two-field-formulation}.

For a lattice $\Lambda=A\mathds Z^{d}$, we will use $\bm x,\bm y\in
\Lambda$ to denote lattice positions and equip the notation with a
tilde to denote the associated integer coordinates $\tilde {\bm x} =
A^{-1}\bm x\in \mathds Z^d$. Further, $\bm p,\bm q\in E^\ast$ are
used for wave vectors, with reduced vectors on the torus $\tilde {\bm
p}=A^{\top}\bm p \in \mathds T^d$. The symbol $A$ is reserved for the
lattice matrix, while $B$ is used for matrices of other dimensions,
in particular for the Galerkin matrix in the following subsections.

We then briefly discuss Fourier functions, in whose basis the
convolution operator is diagonalised, as well as Sobolev spaces.
\begin{definition}[Fourier functions, Sobolev spaces]
  The Fourier functions
  $$
  v_{\bm x}(\bm p)=e^{-\imath 2\pi \bm x\cdot \bm p}\,,\ \bm x\in \Lambda,
  $$
  form an orthonormal basis of $L^2(E^\ast)$ with $\langle v_{\bm
  x},v_{\bm y}\rangle= \delta_{\bm x,\bm y}$ with respect to the
  inner product from Def.~\ref{def:gap_space} where $\delta_{\bm
  x,\bm y}$ denotes the Kronecker delta. Then every $f\in
  L^2(E^\ast)$ has the expansion
  $$
  f(\bm p)=\sum_{\bm x\in \Lambda} v_{\bm x} f_{\bm x}
  $$
  with Fourier coefficients
  \begin{equation}
    \label{eq:fourier_coeffs}
    f_{\bm x}= \langle v_{\bm x},f \rangle=(\mathcal F_\Lambda^{-1}
    f)(\bm x),\quad \bm x\in \Lambda,
  \end{equation}
  for which Parseval's identity
  $$
  \langle f,g \rangle=\sum_{\bm x\in \Lambda} \bar f_{\bm x} g_{\bm x}
  $$
  holds. In the reduced coordinates $\tilde {\bm p}=A^{\top}\bm p \in
  \mathds T^d$, this corresponds to the standard Fourier series on
  the torus. Finally, for $s\in \mathds R$, the Sobolev space
  $H^s=H^s(E^\ast)$ is defined as the space of all functions $f$ on
  $E^\ast$ with finite norm
  \[
    \Vert f\Vert_{H^s}=\Big(\sum_{\bm x\in \Lambda} (1+\vert \bm
    x\vert^2)^s \vert f_{\bm x}\vert^2\Big)^{1/2},
  \]
  thus $H^0=L^2(E^\ast)$. In the case of vector- or matrix-valued
  functions $f$, we replace the absolute value in the formula above
  by the Euclidean or the Frobenius norm.
\end{definition}

Then define the convolution operator from
Eq.~\eqref{eq:two-field-formulation}  as
\begin{equation}
  \label{eq:conv_op}
  \mathcal B[f]=\hat K\ast f,
\end{equation}
where, for matrix-valued $f$, the convolution is applied
component-wise. Following the notation for pseudodifferential
operators, the symbol of the operator is given by $\sigma_{\mathcal
B}(\bm x)=K(\bm x)$, $\bm x\in \Lambda$.
\begin{lemma}
  The Fourier functions $\{v_{\bm x}\}_{\bm x\in \Lambda}$ define the
  eigenfunctions of $\mathcal B$ with
  \[
    \mathcal B[v_{\bm x}]=\lambda_{\bm x} v_{\bm x},\quad
    \lambda_{\bm x}=\sigma_{\mathcal B}(\bm x)=K(\bm x),\quad \bm x\in \Lambda.
  \]
\end{lemma}

The proof follows directly from the definition of the Fourier
coefficients. As an immediate consequence, we have for a kernel $\hat
K$ with bound on the Fourier coefficients $|K(\bm x)|\le c(1+\vert
\bm x\vert)^{-\nu}$, $c>0$, for some $\nu>0$, that $\mathcal B$
increases Sobolev regularity by $\nu$, namely
\[
  \mathcal B:H^s\to H^{s+\nu},\quad s\in \mathds R,
\]
which follows from the definition of the Sobolev norm. For kernels of
the form from Eq.~\eqref{eq:interaction},
\[
  K(\bm x)=a(\bm x)+\sum_{j=1}^\ell b_j \mathcal K_{\nu_j}(\bm x),\quad \nu_j>0
\]
i.e. inverse power-laws with exponents $\nu_j$ plus a compactly
supported function $a$, this holds for $\nu=\min_j \nu_j$. Thus the
gap $\hat \Delta=\hat K\ast G$ is $\nu$ orders smoother than $G$. On
the other hand, $G$ is at most as smooth as $\xi$ and $\hat \Delta$,
and can become discontinuous at nodes $\bm p$ where $\xi$ vanishes
and $\hat \Delta$ loses rank, yet always remains bounded by
Lem.~\ref{lem:aux_lemma_gb} and thus square integrable. Our
quadrature accounts for this fact.

\subsection{Periodic B-Splines}
Periodic B-splines of arbitrary order are ideal for providing a
high-quality numerical approximation of periodic functions. The first
B-splines were introduced by the Russian mathematician
Nikolai Lobachevsky of Kazan University in 1842 as part of his
research into errors in astronomical measurements
\cite{ALLASIAandCAVORETTOandDEROSSI:2013}. The term \emph{B-splines},
however, was coined by Isaac Jacob Schoenberg in 1967
\cite{schoenberg1967spline}, corresponding to their role as basis functions, after he and Curry had introduced them on arbitrary knots \cite{CURRYandSCHOENBERG:1966}.


First, we consider the one-dimensional case, $d=1$, and the following
discretisation of the interval $[0,1)$
\begin{eqnarray}
  \label{Eqn: Discretisation [0,1)}
  [0,1)=\bigcup_{\ell=1}^n\,[p_{\ell},p_{\ell+1})\,,\quad
  p_{\ell}=(\ell-1)h\,,\quad h=1/n\,,\quad n\in\mathbb{N}\,,
\end{eqnarray}
where $n$ and, therefore, $h$ will be the same for all space
dimensions. A possible basis
in the space of piecewise polynomials on $[0,1)$ with respect to the
decomposition
(\ref{Eqn: Discretisation [0,1)}) and periodically extended to the
whole $\mathbb{R}$ is denoted by
\begin{eqnarray}
  \label{Eqn: Basis B-splines}
  \Phi_n^{(\kappa)}=\Big(\phi_1^{(\kappa)},\dots,\phi_n^{(\kappa)}\Big)\,,
\end{eqnarray}
where $\kappa$ denotes the degree of the splines. The functions
$\phi_{\ell}^{(\kappa)}$ will be defined recursively
$$
\phi_1^{(0)}(p)=
\begin{cases}
  1 & \text{for } -h/2+m\leq p < h/2+m\,,\quad m\in\mathbb{Z}\,,\\
  0 & \text{elsewhere} \\
\end{cases}\,.
$$
For $\kappa\geq 1$ and for $n>\kappa+1$, we define the B-splines of
the degree $\kappa$ by the use of the recursion
$$
\phi_1^{(\kappa)}(p)=\frac{1}{h}\,\int\limits_{-1/2}^{1/2}
\phi_1^{(\kappa-1)}(q)\phi_1^{(0)}(p-q)\,\mathrm{d}q\,,\quad \kappa=1,2,\dots\ ,
$$
and
$$
\phi_{\ell}^{(\kappa)}(p)=\phi_1^{(\kappa)}(p-p_{\ell})\,,\quad
\ell=2,\dots,n\,,
$$
as well as by the periodical extension
$$
\phi_{\ell}^{(\kappa)}(p+m)=\phi_{\ell}^{(\kappa)}(p)\,,\quad m\in\mathbb{Z}\,.
$$
The finite dimensional space of all periodic piecewise polynomials on
$\mathbb{R}$
will be denoted by
$$
H_n=\mathrm{span}\,\Phi_n^{(\kappa)}=
\big\{u_n=\Phi_n^{(\kappa)}c\,,\quad c\in\mathbb{R}^n\big\}\,.
$$
The basic properties of the B-splines defined above are summarised in
the following theorem
\begin{theorem}
  \label{Thm: B-splines Properties}
  Let $e=(1,\dots,1)^{\top}\in\mathbb{R}^n$.
  \begin{enumerate}
    \item Decomposition of unity:
      $$
      \Phi_n^{(\kappa)}e=\sum_{\ell=1}^n\phi_{\ell}^{(\kappa)}(p)=1\,,\quad
      \text{for all }
      p\in\mathbb{R}\,,\quad \kappa=0,1,\dots\,,\quad n>\kappa+1\,.
      $$
    \item The first B-splines are even functions:
      $$
      \phi_1^{(\kappa)}(p)=\phi_1^{(\kappa)}(-p)\,,\quad \text{for all}\quad
      p\in\mathbb{R}\,,\quad \kappa=0,1,\dots\,.
      $$
    \item The Fourier series:
      $$
      \phi_{\ell}^{(\kappa)}(p)=\sum_{x\in\mathbb{Z}}c_{\ell}^{(\kappa)}(x)\,e^{-\imath\,2\pi\,
      x\, p}\,,
      $$
      where
      $$
      c_1^{(\kappa)}(x)=h\,\mathrm{sinc}^{\kappa+1}(\pi h\,x)\,,
      $$
      with
      $$
      \mathrm{sinc}(p)=
      \begin{cases}
        1 & \text{for}\quad p=0\,, \\[1ex]
        \dfrac{\sin(p)}{p} & \text{for}\quad p\neq 0\,,\quad p\in\mathbb{R} \\
      \end{cases}
      $$
      and
      $$
      c_{\ell}^{(\kappa)}(x)=c_1^{(\kappa)}(x)e^{\imath\,2\pi
      h\,x(\ell-1)}\,,\quad \ell=2,\dots,n\,.
      $$
    \item Sobolev spaces:
      $$
      \phi_{\ell}^{(\kappa)}\in H^a\,,\quad  \text{for all}\quad a<\kappa+1/2\,.
      $$
    \item Approximation property:
      $$
      \inf_{u_n\in H_n}\|u-u_n\|_{H^a}\leq
      c_1\,h^{s-a}\|u\|_{H^s}\,,\quad \text{for all}\quad u\in
      H^s\,,\quad c_1>0\,,
      $$
      where
      $$
      -\infty < a \leq s\leq \kappa+1\,,\quad a < \kappa+1/2\,.
      $$
    \item Inverse inequality:
      $$
      c_2\,h^{s-a}\|u_n\|_{H^s}\leq \|u_n\|_{H^a}\leq
      \|u_n\|_{H^s}\,,\quad  \text{for all}\quad  u_n\in H^s\,,
      $$
      where
      $$
      -\infty < a \leq s < \kappa+1/2\,.
      $$
  \end{enumerate}
\end{theorem}

Now consider the discretisation parameter $n\in \mathds N$ and let
\begin{equation}
  \label{eq:Q_n}
  Q_n=\big\{\bm \ell \in \mathds N^d,\quad \bm
  \ell_j\in\{1,\dots,n\},\quad j=1,\dots, d\big\},
\end{equation}
denote the set of $N=n^d$ $d$-dimensional indices. From now on, we
use $\bm \ell,\bm m,\bm n\in Q_n$ for such multi-indices.

In the multidimensional case $d>1$, and for $\bm p\in E^\ast$, we simply define
\begin{equation}
  \label{Eqn: Multidimensional spline}
  \phi_{\boldsymbol{\ell}}^{(\kappa)}(\boldsymbol{p})=
  \prod_{j=1}^d\phi_{\bm \ell_j}^{(\kappa)}((A^\top \bm p)_j)\,,\quad
  \boldsymbol{\ell}\in Q_n\,,\quad \boldsymbol{p}\in E^\ast
\end{equation}
following by the obvious periodic continuation for all $\bm
p\in\mathbb{R}^d$. Through the substitution $\bm p=A^{-\top} \tilde
{\bm p}$ in Eq.~\eqref{eq:fourier_coeffs}, their Fourier coefficients
factorise into the one-dimensional ones, and we obtain for $\bm x\in
\Lambda$ and $\tilde{\bm x}=A^{-1}\bm x\in \mathbb{Z}^d$ that
\begin{equation}
  \label{eq:fourier_coeffs_spline}
  c_{\boldsymbol{\ell}}^{(\kappa)}(\boldsymbol{x})=\prod_{j=1}^dc_{\bm
  \ell_j}^{(\kappa)}(\tilde{\bm x}_j)=c_{\bm e}^{(\kappa)}(\bm
  x)e^{\imath 2\pi h \tilde {\bm x}\cdot(\bm \ell-\bm e)}\,,\quad
  \boldsymbol{\ell}\in Q_n\,,
\end{equation}
with $\bm e=(1,\dots,1)^\top\in Q_n$ and $c_{\bm e}^{(\kappa)}(A
\tilde {\bm x})=h^d \prod_{j=1}^d \mathrm{sinc}^{\kappa+1}(\pi h
\tilde {\bm x}_j)$.

Thus, the generalised multi-dimensional B-splines can be used for the
high-quality numerical approximation of
periodic pseudodifferential equations. 

\subsection{Galerkin Discretisation}

Let $\kappa=0,1,\dots$ be the degree of the B-splines, and
$n\in\mathbb{N}$ be a one-dimensional discretisation parameter such
that $n>\kappa+1$.
We consider the $n^d-$dimensional B-spline subspace
$$
H_n=\mathrm{span}\,\Phi_n^{(\kappa)}=
\mathrm{span}\,\Big\{\phi_{\bm\ell}^{(\kappa)}\Big\}_{\bm\ell\in Q_n}\,,
$$
with the spline functions $\phi_{\bm\ell}^{(\kappa)}$ defined in
\eqref{Eqn: Multidimensional spline} and $Q_n$ from \eqref{eq:Q_n}.
Note that the knots and therefore also $H_n$ are invariant under the
reflection $\bm p\mapsto -\bm p$. $X_n=H_n^{2\times 2}\cap X$ is then
the space of discrete gap functions with the symmetry described in
Sec.~\ref{sec:model}.

The Galerkin formulation of the system \eqref{eq:two-field-formulation}
reads as follows.
Find $\hat \Delta_n\in X_n$ and $G_n\in X_n$ such that
\begin{eqnarray*}
  \langle v, \hat \Delta_n \rangle &=&\langle v,\mathcal{B}[G_n]\rangle \,, \\
  \langle v, G_n \rangle &=& \langle v, G_\beta[\hat \Delta_n]\rangle
  \quad\text{for all}\quad v\in H_n\,.
\end{eqnarray*}
where the product $\langle \cdot,\cdot\rangle$ with a matrix argument
is applied component-wise. Since $G_\beta$ maps $X$ to $X$ by
Lem.~\ref{lem:aux_lemma_gb}, and both the convolution with the even
kernel and the $L^2$-projection on $H_n$ preserve the symmetry, the
right-hand sides all lie in $X_n$ if $\hat \Delta_n$ and $G_n$ do. We
therefore impose the correct symmetry at the initial guess only, and
it is preserved under iteration.
With the coefficient matrices $\underline \Delta_{\bm m},\underline
G_{\bm m}\in \mathds C^{2\times 2}$
$$
\hat \Delta_n=\Phi_n^{(\kappa)}\underline{\Delta}=
\sum_{\bm m\in Q_n}\underline{\Delta}_{\,\bm m}\,\phi_{\bm m}^{(\kappa)}\,,\quad
G_n=\Phi_n^{(\kappa)}\underline{G}=
\sum_{\bm m\in Q_n}\underline{G}_{\,\bm m}\,\phi_{\bm m}^{(\kappa)},
$$
viewed as columns of length $N$, the Galerkin formulation is given by
the following nonlinear algebraic system
\begin{equation}
  \label{Eqn: Nonlinear algebraic system}
  M\underline{\Delta}=B\underline{G}\,,\quad
  M\underline{G}=\underline g(\underline \Delta),\quad \underline
  g_{\bm \ell}(\underline \Delta)=\langle \phi_{\bm
  \ell}^{(\kappa)},G_\beta[\hat \Delta_n] \rangle,\quad \bm \ell \in Q_n,
\end{equation}
with the mass matrix $M\in\mathbb{R}^{N\times N}$ and the Galerkin matrix
$B\in\mathbb{R}^{N\times N}$ given as
$$
M_{\bm \ell,\bm m}=\langle \phi_{\bm \ell}^{(\kappa)},\phi_{\bm
m}^{(\kappa)}\rangle \,,\quad
B_{\bm \ell,\bm m}=\langle\phi_{\bm m}^{(\kappa)}, \mathcal{B}\big[\phi_{\bm
\ell}^{(\kappa)}\big]\rangle\,,\quad
\bm \ell,\bm m\in Q_n\,.
$$
By Parseval´s theorem, their entries are lattice sums over Fourier coefficients
\begin{equation}
  \label{eq:matrix_coeffs_galerkin}
  M_{\bm \ell,\bm m}=\sum_{\bm x\in \Lambda} \overline{c_{\bm
  \ell}^{(\kappa)}(\bm x)} c_{\bm m }^{(\kappa)}(\bm x),\qquad B_{\bm
  \ell,\bm m}=\sum_{\bm x\in \Lambda} \overline{c_{\bm
  \ell}^{(\kappa)}(\bm x)}\sigma_{\mathcal B}(\bm x) c_{\bm
  m}^{(\kappa)}(\bm x).
\end{equation}

The mass matrix is symmetric, positive definite, block-circulant and
extremely sparse. The Galerkin matrix is real symmetric, since both
$K$ and $|c_{\bm e}^{(\kappa)}|^2$ are real and even,
block-circulant, and positive-semidefinite for $K\ge 0$ and positive
definite for $K>0$. The right-hand side $\underline g(\underline
\Delta)$ needs to be computed through quadrature. For the
mass-matrix, Gauss-quadrature with $\kappa+1$ points per direction is
exact, whereas the nonlinearity is integrated to an accuracy
depending on its regularity. This can be improved using subdivision
in the vicinity of the gap nodes.

The nonlinear algebraic system \eqref{Eqn: Nonlinear algebraic
system} can be solved in different iterative ways. All classical
methods, from simple fixed-point iteration to quasi-Newton or full Newton
methods, can be applied, while taking precautions, based on the
energy functional, to prevent convergence to saddle points or local
maxima. As an example, the fixed-point iteration
\[
  M \underline \Delta_{j+1}=B \underline G_j,\quad M \underline
  G_j=\underline g(\underline \Delta_j),
\]
is the Galerkin form of the self-consistent iteration described in
Prop.~\ref{prop:convergence}, and never increases the discrete free
energy while driving the discrete residual to zero. Convergence
acceleration, i.e. through Broyden updates, is possible. Note that
every step in the solution requires products with $B$ and solutions
with $M$, with both matrices block-circulant, as shown in the next subsection.

\subsection{Block-Circulant Matrices and Eigenvalues of the Galerkin Matrices}

We call a matrix $B\in \mathds C^{N\times N}$, $N=n^d$, indexed by
$Q_n$, block-circulant, if $B_{\bm \ell \bm m}=b(\bm \ell-\bm m)$ for
an $n$-periodic function $b:\mathds Z^d\to \mathds C$. In case that
$d=1$, they are the classical circulant matrices, see
\cite{RJASANOW17} for a detailed description of their properties and
for effective algorithms. These matrices are diagonalised by the
discrete Fourier transform.

\begin{lemma}
  Consider a $n$-periodic function $b\,:\,\mathbb{Z}^d\rightarrow\mathbb{C}$
  $$
  b(\bm \ell+n\bm m)=b(\bm \ell)\quad
  \text{for all}\quad \bm \ell,\bm m\in\mathbb{Z}^d
  $$
  and a block-circulant matrix
  $$
  B\in\mathbb{C}^{N\times N}\,,\quad N=n^d
  $$
  with
  $$
  B_{\bm \ell\bm m}=b(\bm \ell-\bm m)\,,\quad
  \bm \ell,\bm m\in Q_n\,.
  $$
  Furthermore, consider the $d$-dimensional Fourier vectors
  $$
  v_{\bm \ell}\in\mathbb{C}^{N}\,,\quad \bm \ell\in Q_n\,,\quad N=n^d
  $$
  with the components
  $$
  \big(v_{\bm \ell}\big)_{\bm m}=e^{\imath\, 2\pi h\,
  (\bm \ell-\bm e)\cdot(\bm m - \bm e)}\,,\quad
  \bm \ell,\bm m\in Q_n\,,\quad \bm e=(1,\dots,1)^{\top}\in Q_n\,,\quad
  h=\dfrac{1}{n}\,.
  $$
  The vectors $ v_{\bm \ell}$ are the eigenvectors of the matrix $B$
  and the corresponding eigenvalues $\lambda_{\bm \ell}$ are
  $$
  \lambda_{\bm \ell}=
  \sum_{\bm m\in Q_n}B_{\bm e\bm m}\big(v_{\bm \ell}\big)_{\bm m}\,,\quad
  \bm \ell\in Q_n\,.
  $$
\end{lemma}
\begin{proof}\\
  We compute
  \begin{eqnarray*}
    \big(Bv_{\bm \ell}\big)_{\bm k}&=&\sum_{\bm m\in Q_n}B_{\bm k\bm m}
    \big(v_{\bm \ell}\big)_{\bm m}=\sum_{\bm m\in Q_n}b(\bm k-\bm m)
    \big(v_{\bm \ell}\big)_{\bm m}\\[1ex]
    &=&\sum_{\hat{\bm m}\in Q_n-\bm k+\bm e}b(\bm e-\hat{\bm m})
    \big(v_{\bm \ell}\big)_{\hat{\bm m}+\bm k-\bm e}
    =\big(v_{\bm \ell}\big)_{\bm k}\sum_{\bm m\in Q_n}B_{\bm e\bm m}
    \big(v_{\bm \ell}\big)_{\bm m}\,.
  \end{eqnarray*}
  Here, the substitution $\bm \hat{\bm m}=\bm m-\bm k+\bm e$ and the property
  $$
  \big(v_{\bm \ell}\big)_{\bm m+\bm k-\bm e}=
  \big(v_{\bm \ell}\big)_{\bm m}\,\big(v_{\bm \ell}\big)_{\bm k}\,,
  $$
  have been used, and the hat-sign has been removed. Using the
  $n$-periodical property of the function $b$ and of the components
  $\big(v_{\bm \ell}\big)_{\bm m}$, the summation domain was shifted
  back to the original one.
  Thus, we get
  $$
  Bv_{\bm \ell}=\lambda_{\bm \ell}v_{\bm \ell}\,,\quad \bm \ell\in Q_n\,.
  $$
\end{proof}

Thus, a product with a block-circulant matrix and the solution of a
system with a regular block-circulant matrix require only one
$d$-dimensional Fast Fourier transform pair of length $N$ each with
log-linear cost scaling in $N$.

\begin{proposition}
  \label{Ex: Operator Galerkin B-Splines Convolution}
  Both the mass matrix $M$ and the Galerkin matrix $B$ are block-circulant, with
  $$
  B_{\boldsymbol{\ell}\boldsymbol{m}}=b(\boldsymbol{\ell}-\boldsymbol{m})=
  \sum_{\boldsymbol{x}\in\Lambda}
  \left|c_{\boldsymbol{e}}^{(\kappa)}(\boldsymbol{x})\right|^2
  \sigma_{\mathcal{B}}(\boldsymbol{x})
  e^{\imath\,2\pi h\,\tilde{\bm x}\cdot(\bm \ell-\bm m)}\,,\quad
  \bm \ell,\bm m\in Q_n\,,\quad \tilde{\bm x}=A^{-1}\bm x\,.
  $$
  and for $M$ likewise with $\sigma_{\mathcal B}$ replaced by 1. The
  eigenvectors are given by the Fourier vectors $\bm v_{\ell}$, and
  the eigenvalues of $B$ read for $\bm \ell\in Q_n$,
  \begin{equation}
    \label{eq:final_eq_eigenvalues}
    \lambda_{\bm \ell}=h^d \sum_{\bm x\in \mathds Z^d}\prod_{j=1}^d
    \mathrm{sinc}^{2(\kappa+1)}\Big(\pi h(\bm \ell_j-1+n \bm
    x_j)\Big)\sigma_{\mathcal B}(A(\bm \ell-\bm e+n\bm x)),
  \end{equation}
  and $\bm e=(1,\dots,1)^\top$.
  For indices such that $\bm \ell_j\neq 1$ in every direction, this yields
  \[
    \lambda_{\bm \ell}\!=\!
    h^d\!\prod_{j=1}^d\!\!\left(\!\dfrac{\sin(\pi h\,(\bm
    \ell_j\!-\!1))}{\pi h}\!\right)^{\!2(\kappa+1)}
    \!\!\!\sum_{\bm x \in\mathbb{Z}^d}\!
    \!\left(\prod_{j=1}^d\!\dfrac{1}{|\bm \ell_j\!-\!1\!+\!n\bm
    x_j|^{2(\kappa+\!1)}}\!\right)\!\!
    \sigma_{\mathcal{B}}\big(A(\bm \ell-\bm e+n\bm x)\big).
  \]
  In case that $\bm \ell_j=1$, the sum over the coordinate $\bm x_j$
  only contributes the term $\bm x_j=0$. In particular, $\lambda_{\bm
  e}=h^d K(\bm 0)$. The eigenvalues of the mass-matrix $M$ are given
  by Eq.~\eqref{eq:final_eq_eigenvalues} under the replacement
  $\sigma_{\mathcal B}=1$ and factorise into the product of the
  one-dimensional eigenvalues.
\end{proposition}
\begin{proof}
  We begin with the entries of the matrix $B$ from
  Eq.~\eqref{eq:matrix_coeffs_galerkin} and apply the shift formula
  for the spline Fourier coefficients in
  Eq.~\eqref{eq:fourier_coeffs_spline}. With $\tilde {\bm x}=A^{-1}
  \bm x$, we obtain
  $$
  B_{\bm e \bm m}=
  \sum_{\bm x \in\Lambda}
  \left|c_{\boldsymbol{e}}^{(\kappa)}(\bm x)\right|^2
  \sigma_{\mathcal{B}}(\bm x)
  e^{\imath\,2\pi h\,\tilde{\bm x}\cdot(\bm e-\bm m)}\,.
  $$
  Therefore, the eigenvalues are
  $$
  \lambda_{\bm \ell}=
  \sum_{\bm x \in\Lambda}
  \left|c_{\boldsymbol{e}}^{(\kappa)}(\bm x)\right|^2
  \sigma_{\mathcal{B}}(\bm x)
  \sum_{\bm m\in Q_n}
  e^{\imath\, 2\pi h\, (\bm m-\bm e)\cdot(\bm \ell-\bm e-\tilde{\bm x})}\,.
  $$
  Since
  $$
  \sum_{\bm m\in Q_n}
  e^{\imath\,2\pi h\,(\bm m-\bm e)\cdot(\bm \ell-\bm e-\tilde{\bm x})}=
  \begin{cases}
    n^d &\quad\text{for}\quad\tilde{\bm x}=\bm \ell-\bm e+n\hat{\bm
    x}\,,\quad  \hat{\bm x}\in\mathbb{Z}^d\\[1ex]
    0 & \quad\text{else} \\
  \end{cases}\,,
  $$
  and removing the hat-sign, we get
  $$
  \lambda_{\bm \ell}=
  n^d\sum_{\bm x \in\mathbb{Z}^d}
  \left|c_{\boldsymbol{e}}^{(\kappa)}\big(A(\bm \ell-\bm e+n\bm x)\big)\right|^2
  \sigma_{\mathcal{B}}\big(A(\bm \ell-\bm e+n\bm x)\big)\,.
  $$
  Remember
  $$
  c_{\boldsymbol{e}}^{(\kappa)}\big(A(\bm \ell-\bm e+n\bm x)\big)=
  h^d\prod_{j=1}^d\mathrm{sinc}^{\kappa+1}(\pi h\,(\bm \ell_j-1+n\,\bm x_j))\,.
  $$
  This yields Eq.~\eqref{eq:final_eq_eigenvalues}, and the following
  equation is obtained from
  \[
    \sin\big(\pi h(\bm \ell_j-1+n\bm x_j)\big)=\pm \sin\big(\pi h(\bm
    \ell_j-1)\big),
  \]
  whereas for $\bm \ell_j=1$, we have
  $\mathrm{sinc}^{2(\kappa+1)}(\pi \bm x_j)=\delta_{\bm x_j,0}$.
  Finally, for $\sigma_{\mathcal B}=1$, each summand and thus also
  the sum is a product over $j$.
\end{proof}

\begin{example}
  \label{Ex: Eigenvalues Galerkin Zeta}
  Consider the integer lattice $\Lambda=\mathds Z$ and the kernel
  $\hat K=C_0+U_0 Z_{\Lambda,\nu}$, which will be used in analytical
  solutions to the gap equation in the next section. The symbol then reads
  $$
  \sigma_\mathcal{B}(x)=\dfrac{U_0}{|x|^\nu}\,,\quad \nu>0\,,\quad
  x\neq 0,\quad \sigma_{\mathcal B}(0)=C_0.
  $$
  Then Prop.~\ref{Ex: Operator Galerkin B-Splines Convolution} yields
  $\lambda_1=h C_0$ and the remaining eigenvalues
  $\lambda_\ell\,,\ell>1$ are of the form
  $$
  \lambda_\ell=h U_0\left(\dfrac{\sin(\pi h\,(\ell-1))}{\pi
  h}\right)^{2(\kappa+1)}
  \sum_{x\in\mathbb{Z}}\dfrac{1}{|\ell-1+xn|^{2(\kappa+1)+\nu}}\,.
  $$
  This series converges reasonably quickly. However, it can be
  expressed through the Epstein zeta function as follows, using the
  following abbreviation $\eta=2(\kappa+1)+\nu$
  $$
  \lambda_\ell=h U_0\left(\dfrac{\sin(\pi h\,(\ell-1))}{\pi
  h}\right)^{2(\kappa+1)}
  Z_{n\mathbb{Z},\eta}(1-\ell, 0).
  $$
  These eigenvalues can be efficiently precomputed with the
  help of EpsteinLib \cite{epsteinlib}.
\end{example}

In dimensions $d\ge 2$, the lattice sum forming the eigenvalues is no
longer an Epstein zeta function. However, the summand decays as
$\vert \bm x_j\vert^{-2(\kappa+1)-\nu}$, hence the sum can be
evaluated directly. Importantly, the eigenvalues can be precomputed,
and can then be reused at every iteration step.

With the eigenvalues of both $M$ and $B$ given, a single iteration
step of a fixed-point iteration requires the quadrature of the
nonlinearity, thus $\mathcal O(N)$ evaluations of $G_\beta$, and a
small number of fast Fourier transforms. Meanwhile, memory only grows
linearly with $N$. We report on the measured numerical work and
memory of our algorithm in detail in the next section.

\section{Non-Trivial Analytical Solutions from Band Engineering}
\label{sec:analytical}

In this section, we derive a parametrised family of solutions to the
gap equation with long-range interactions for the one-dimensional
lattice $\Lambda=\mathds Z$ with an attractive interaction $K\ge 0$,
whose lattice Fourier series reads
\[
  \hat K(p) = C_0+U_0 Z_{\mathds Z,\nu}(p),
\]
where $C_0,U_0\ge 0$ and $\nu>0$.
This is achieved through band-engineering, meaning that the
dispersion relation $\xi$
is fine-tuned such that a specified trial gap function $\Delta$
becomes a solution to the gap equation.
While trivial solutions, such as constant solutions for
nearest-neighbor interactions,
are known in the literature, no closed-form solutions exist for
singular power-law interactions to date. In this section, we fill
this gap, providing a reliable benchmark for solvers in
unconventional superconductivity.

We begin by presenting an even solution, which is referred to as a
singlet solution.
\begin{lemma}[1D singlet solution]
  \label{lem:singlet1d}
  Let $\Lambda=\mathds Z$, $\beta=\infty$, $C_0>0$, $U_0\ge 0$,
  $\nu>0$, and $\eta>1$. Set the anomalous density $\alpha:E^\ast \to\mathds R$,
  \[
    \alpha(p)=\alpha_0+\frac{\alpha_1}{2}\big(Z_{\Lambda,\eta}(p-1/4)+Z_{\Lambda,\eta}(p+1/4)\big),
  \]
  with
  \[
    \frac{-Z_{\Lambda,\eta}(\tfrac{1}{4}
    )}{Z_{\Lambda,\eta}(0)+Z_{\Lambda,\eta}(\tfrac{1}{2})-2Z_{\Lambda,\eta}(\tfrac{1}{4})}<\alpha_0<\frac{1}{2},\quad
    \alpha_1=\frac{1-2\alpha_0}{Z_{\Lambda,\eta}(0)+Z_{\Lambda,\eta}(\tfrac{1}{2})}.
  \]
  Further fix $C_0, U_0$ such that
  \[
    C_0\alpha_0>-U_0 \alpha_1 Z_{\Lambda,\nu+\eta}(\tfrac{1}{4}),
  \]
  and set the singlet gap $\psi:E^\ast\to \mathds R$ and the
  dispersion relation $\xi:E^\ast \to \mathds R$ as
  \begin{align*}
    \psi(p) &= C_0 \alpha_0+\frac{U_0
    \alpha_1}{2}\big(Z_{\Lambda,\eta+\nu}(p-1/4)+Z_{\Lambda,\eta+\nu}(p+1/4)\big),\\
    \xi(p)&=\mathrm{sgn}(|p|-1/4)\psi(p)
    \frac{\sqrt{1-4\alpha(p)^2}}{2\alpha(p)}.
  \end{align*}
  Then, $0<\alpha(p)\le 1/2$ with equality at $p=\pm 1/4$, $\xi$ is
  even and continuous with $\xi(\pm \tfrac{1}{4})=0$, and, finally,
  \[
    \hat \Delta(p)=
    \begin{pmatrix}
      0 & \psi(p)\\
      -\psi(p) & 0
    \end{pmatrix}\in L^2(E^\ast,\mathds C^{2\times 2})
  \]
  is a solution to the gap equation from Def.~\ref{def:gap_eq} with
  kernel $\hat K(p)=C_0+U_0Z_{\Lambda,\nu}(p)$ and degenerate
  eigenvalues of $E[\hat \Delta]$ equal to $\lambda=\psi/(2\alpha)$.
\end{lemma}

The proof of above lemma requires several technical lemmas.

\begin{lemma}\label{lem:zeta-conv}
  The family of zeta functions $(Z_\nu)_{\nu \in \mathds C}$ with
  $\Lambda = \mathds Z$ forms a group under periodic convolution with
  \[
    Z_{\nu_1} \ast Z_{\nu_2} = Z_{\nu_1 + \nu_2},
  \]
  where the convolution is defined in the sense of distributions in
  case $\nu_1$ and $\nu_2$ have both no real part larger than unity.
\end{lemma}

\begin{proof}
  For $\nu_1, \nu_2$ with real part larger than unity this follows
  directly as the zeta functions are lattice Fourier transforms of
  the homogeneous regularised power-law $\mathcal K_\nu$, for which
  $\mathcal K_{\nu_1} \mathcal K_{\nu_2}=\mathcal K_{\nu_1+\nu_2}$ holds.
  To cover the whole complex plane we note that both sides are weakly
  holomorphic and hence the identity extends to $\mathds C \times \mathds C$.
\end{proof}

Next we turn to a closer analysis of the sign and the monotonicity of
the zeta function in one dimension.

\begin{lemma}\label{lem:zeta-monotone}
  For real $\nu > 1$, $Z_\nu$ is even, continuous and strictly
  decreasing on $[0, 1/2]$.
\end{lemma}

\begin{proof}
  Owing to $\nu > 1$, we can write the series for $Z_\nu$ as a series
  over Mellin transform of the decaying exponential,
  \begin{equation}\label{eq:zeta-poisson}
    Z_\nu(p) = 2 \sum_{m = 1}^\infty \frac{\cos(2 \pi m p)}{m^\nu}
    = \frac{1}{\Gamma(\nu)} \int_0^\infty t^{\nu - 1} \big(
      P_{e^{-t}}(2 \pi p) - 1
    \big)
    \, \text d t,
  \end{equation}
  with the Poisson kernel
  \[
    P_r(\theta) = \sum_{m = 0}^\infty r^m \cos \theta = \frac{1 -
    r^2}{1 - 2 r \cos \theta + r^2},
    \quad r \in (0, 1), \quad \theta \in \mathds R.
  \]
  Since the kernel $P_r$ is strictly decreasing in this argument on
  $[0, \pi]$ for any $0 < r < 1$, it follows that
  $Z_\nu$ as parameter integral over the Poisson kernel is a
  decreasing function on $[0, 1/2]$.
  For $0 \leq p_1 < p_2 \leq 1/2$, the difference $P_r(2 \pi p_1) -
  P_r(2 \pi p_2)$ is thus strictly positive
  for any $0 < r < 1$, so also the integral against a positive weight
  like in~\eqref{eq:zeta-poisson} is strictly positive.
  This proves that $Z_\nu$ is strictly decreasing on $[0, 1/2]$.
\end{proof}

\begin{lemma}\label{lem:zeta-sum}
  With
  \[
    S_\nu(p) = \frac{1}{2} \big( Z_\nu(p - 1/4) + Z_\nu(p + 1/4)
    \big), \quad p \in (-1/2, 1/2),
  \]
  it holds $S_\nu(p) = 2^{-\nu} Z_\nu(2 p + 1/2)$.
\end{lemma}

\begin{proof}
  For a parameter $\nu$ with $\mathrm{Re}\,\nu > 1$ this follows
  directly by algrebaic manipulations of the series terms,
  particularly that $\cos(\pi / 2\, m)$ vanishes for odd $m$. The
  identity then extends to the complex plane
  since both sides are weakly holomorphic.
\end{proof}

We now turn to the proof of Lemma~\ref{lem:singlet1d}.

\begin{proof}
  First, we note that the identity for $\psi$ follows immediately from the
  group property of Lemma~\ref{lem:zeta-conv} that simply shifts the Epstein
  zeta parameter. Next, we show that the constraints on $\alpha_0$ and
  $\alpha_1$ guarantee $0<\alpha\le 1/2$. We abbreviate
  $Z=Z_{\Lambda,\eta}$. By Lemma~\ref{lem:zeta-sum},
  \[
    \alpha(p) = \alpha_0 + \alpha_1 S_\eta(p) = \alpha_0 + 2^{-\eta}
    \alpha_1 Z(2 p + 1/2).
  \]
  Evaluating Lemma~\ref{lem:zeta-sum} at $p=1/4$ and $p=0$ gives
  \[
    Z(0)+Z(\tfrac12) = 2^{1-\eta} Z(0) > 0, \qquad
    Z(\tfrac14) = 2^{-\eta} Z(\tfrac12) < 0,
  \]
  since 
  $$
   Z(0)=2\zeta(\eta)>0\quad \text{and}\quad Z(\tfrac12)=-2\sum_{m\ge1}(-1)^{m+1}m^{-\eta}<0.
  $$
  In particular, $\alpha_1>0$ as $\alpha_0<1/2$. As $Z$ is $1$-periodic,
  even, and, by Lemma~\ref{lem:zeta-monotone}, strictly decreasing on
  $[0,1/2]$, it attains its maximum exactly on $\mathds Z$ and its
  minimum exactly on $1/2+\mathds Z$. Hence, for $p\in E^\ast$,
  $\alpha$ is maximal exactly at $p=\pm 1/4$ with
  \[
    \alpha(\pm\tfrac14) = \alpha_0 +
    \frac{\alpha_1}{2}\big(Z(0)+Z(\tfrac12)\big) = \frac12,
  \]
  and minimal at $p=0$ with
  \[
    \alpha(0) = \alpha_0 + \alpha_1 Z(\tfrac14)
    =
    \frac{\alpha_0\big(Z(0)+Z(\tfrac12)-2Z(\tfrac14)\big)+Z(\tfrac14)}{Z(0)+Z(\tfrac12)}
    > 0,
  \]
  where positivity is equivalent to the lower bound on $\alpha_0$.
  Thus $0<\alpha(p)\le 1/2$ with equality only at $p=\pm 1/4$.

  Next, we turn to the dispersion relation $\xi$. Since
  $Z_{\Lambda,s}$ is even and, for $s>1$, continuous, $\alpha$ and
  $\psi$ are even and continuous, and so is $\xi$ away from
  $p=\pm 1/4$, where the sign factor jumps. There, however,
  $\alpha(\pm\tfrac14)=1/2$, so $\sqrt{1-4\alpha(p)^2}\to 0$ as
  $p\to\pm\tfrac14$, while $\psi$ and $1/(2\alpha)\le 1/(2\alpha(0))$
  remain bounded. Hence $\xi$ is even and continuous with
  $\xi(\pm\tfrac14)=0$.

  Lastly, we show that the gap $\psi$ is positive. By Lemma~\ref{lem:zeta-sum}
  with $\eta+\nu$ in place of $\eta$,
  \[
    \psi(p) = C_0\alpha_0 + 2^{-(\eta+\nu)} U_0 \alpha_1
    Z_{\Lambda,\eta+\nu}(2p+1/2).
  \]
  As $U_0\ge 0$ and $\alpha_1>0$, the same monotonicity argument as for
  $\alpha$ shows that $\psi$ is minimal at $p=0$, where
  \[
    \psi(0) = C_0\alpha_0 + U_0\alpha_1 Z_{\Lambda,\eta+\nu}(\tfrac14) > 0
  \]
  by the assumption on the constants $C_0$ and $U_0$. In particular, $\xi$ has the
  sign of $|p|-1/4$, and $\lambda=\psi/(2\alpha)>0$.

  Now that the bounds and the symmetries of the involved function are proved,
  we show that $\psi$ is indeed the singlet solution to the gap equation
  with dispersion $\xi$. To that end, we need to prove that
  \[
    \alpha J = G_\infty[\psi J],
  \]
  with the matrix
  \[
    J =
    \begin{pmatrix}
      0 & 1 \\
      -1 & 0
    \end{pmatrix}.
  \]
  Direct calculation shows
  \[
    G_\infty[\psi J] = \frac{\psi}{2 \sqrt{\xi^2 + \psi^2}} J,
  \]
  so that the equation for $\alpha$ reduces to the identity,
  \[
    4 \alpha^2 (\xi^2 + \psi^2) = \psi^2,
  \]
  that is
  \[
    \xi^2 = \psi^2 \frac{1 - 4 \alpha^2}{4 \alpha^2},
  \]
  precisely following from the definition of the dispersion.
\end{proof}

\begin{remark}[1D triplet solution]
  \label{rem:triplet_1d}
  The construction from Lem.~\ref{lem:singlet1d} has an odd symmetry
  equivalent. Let $\Lambda=\mathds Z$, $\beta=\infty$, with $\eta>1$
  and $\nu>0$, and consider an interaction $\hat K=C_0+U_0
  Z_{\Lambda,\nu}$ with $U_0>0$ and $C_0\in \mathds R$. Define
  \[
    \alpha(p) =
    \frac{Z_{\Lambda,\eta}(p-\tfrac{1}{4})-Z_{\Lambda,\eta}(p+\tfrac{1}{4})}{2\big(Z_{\Lambda,\eta}(0)-Z_{\Lambda,\eta}(\tfrac{1}{2})\big)},\quad
    d(p)=U_0\frac{Z_{\Lambda,\eta+\nu}(p-\tfrac{1}{4})-Z_{\Lambda,\eta+\nu}(p+\tfrac{1}{4})}{2\big(Z_{\Lambda,\eta}(0)-Z_{\Lambda,\eta}(\tfrac{1}{2})\big)},
  \]
  with $d=\hat K\ast \alpha$. Set
  \[
    \xi(p)= \mathrm{sgn}(\vert
    p\vert-\tfrac{1}{4})\frac{d(p)}{2\alpha(p)}\sqrt{1-4\alpha(p)^2},
  \]
  continuously extended to $p\in\{0,\pm\tfrac{1}{2}\}$ of $a$ and
  $d$, with $\xi$ even, continuous, and vanishing at $\pm
  \tfrac{1}{4}$. Then for any $M\in \mathds R^{2\times 2}$ with
  $M^\top=M$ and $M^2=\mathds 1_2$,
  \[
    \hat \Delta(p)=M d(p)\in X,
  \]
  solves the gap equation.
\end{remark}

\pgfplotsset{
  lemmaaxis/.style={
    width=6.4cm, height=4.8cm,
    xmin=-0.5, xmax=0.5,
    xtick={-0.5,-0.25,0,0.25,0.5},
    xticklabels={$-\frac12$,$-\frac14$,$0$,$\frac14$,$\frac12$},
    xlabel={$p$},
    label style={font=\footnotesize},
    tick label style={font=\scriptsize},
    legend style={font=\scriptsize, draw=none, fill=none},
    legend cell align=left,
    yticklabel style={/pgf/number format/fixed, /pgf/number format/precision=2},
  },
  swave/.style={blue!70!black, thick},
  pwave/.style={red!70!black, thick},
}
\newcommand{\LemmaFermi}[2]{%
  \draw[gray, dotted] (axis cs:-0.25,#1) -- (axis cs:-0.25,#2);
  \draw[gray, dotted] (axis cs:0.25,#1) -- (axis cs:0.25,#2);
  \draw[gray, thin] (axis cs:-0.5,0) -- (axis cs:0.5,0);}
\newcommand{\LemmaBandOf}[4]{%
\begin{tikzpicture}
  \begin{axis}[lemmaaxis, ymin=#3, ymax=#4, legend pos=south east]
    \LemmaFermi{#3}{#4}
    \addplot[#1] table[x=p, y=xi] {Figures/#2.dat};
    \addlegendentry{$\xi$}
    \addplot[#1, thin, dashed] table[x=p, y=lambda] {Figures/#2.dat};
    \addlegendentry{$\lambda$}
  \end{axis}
\end{tikzpicture}}
\newcommand{\LemmaBand}{\LemmaBandOf{swave}{Lemma61Solution}{-1}{1.3}}
\newcommand{\LemmaBandP}{\LemmaBandOf{pwave}{Lemma61Triplet}{-0.35}{0.35}}
\newcommand{\LemmaGap}{%
\begin{tikzpicture}
  \begin{axis}[lemmaaxis, ymin=0.69, ymax=0.74,
    ytick={0.69,0.7,0.71,0.72,0.73,0.74}, ylabel={$\psi$}]
    \draw[gray, dotted] (axis cs:-0.25,0.69) -- (axis cs:-0.25,0.74);
    \draw[gray, dotted] (axis cs:0.25,0.69) -- (axis cs:0.25,0.74);
    \addplot[swave] table[x=p, y=psi] {Figures/Lemma61Solution.dat};
  \end{axis}
\end{tikzpicture}}
\newcommand{\LemmaGapP}{%
\begin{tikzpicture}
  \begin{axis}[lemmaaxis, ymin=-0.16, ymax=0.16, ylabel={$d$}]
    \LemmaFermi{-0.16}{0.16}
    \addplot[pwave] table[x=p, y=d] {Figures/Lemma61Triplet.dat};
  \end{axis}
\end{tikzpicture}}

\begin{figure}[ht]
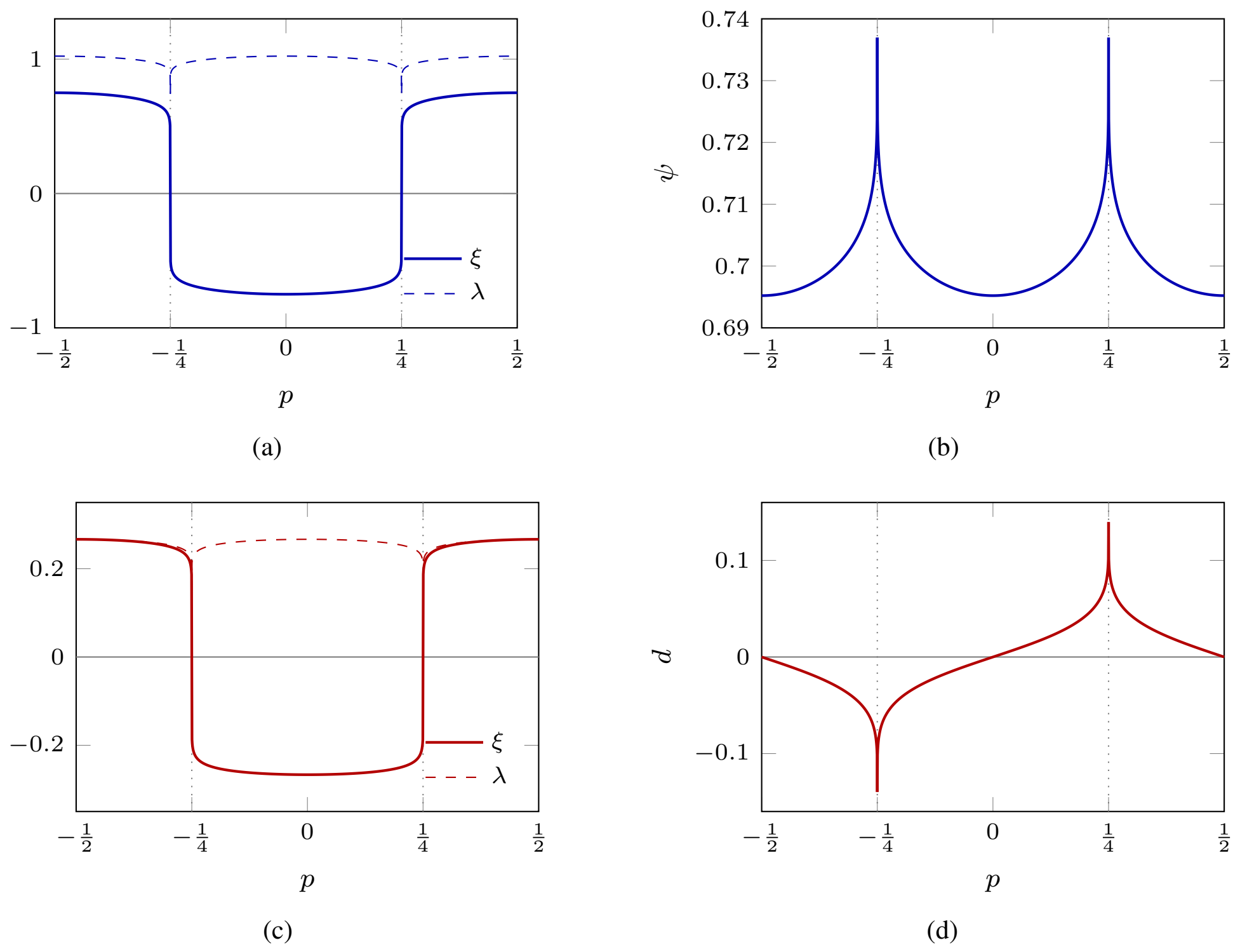

  \centering
  \subfloat[]{\LemmaBand}\hfill
  \subfloat[]{\LemmaGap}\\
  \subfloat[]{\LemmaBandP}\hfill
  \subfloat[]{\LemmaGapP}
  \caption{Exact solutions over $E^\ast = (-1/2, 1/2)$,
    $d=1$, $T=0$, $\eta=1.1$, $\nu=0.1$, $C_0=2$, $U_0=0.5$; dotted:
    Fermi points $p=\pm 1/4$.
    (a), (b) Singlet ($s$-wave) solution of Lemma~\ref{lem:singlet1d},
    $\alpha_0=0.35$: engineered band $\xi$, quasi-particle energy
    $\lambda=\psi/(2\alpha)$, and gap $\psi$.
    (c), (d) Triplet ($p$-wave) solution of
    Remark~\ref{rem:triplet_1d}: band $\xi$, $\lambda=d/(2\alpha)$, and
  gap amplitude $d$.}
  \label{fig:lemmasolution}
\end{figure}

The proof proceeds in analogy to the singlet case in Lem.~\ref{lem:singlet1d}.

Figure~\ref{fig:lemmasolution} shows the singlet solution of
Lemma~\ref{lem:singlet1d} with $\eta=1.1$, $\nu=0.1$, $C_0=2$,
$U_0=0.5$ and $\alpha_0=0.35$: the band $\xi$ changes sign at the
Fermi points $p=\pm 1/4$, where the quasi-particle energy
$\lambda=\psi/(2\alpha)$ is minimal and the gap $\psi$ is maximal.
Near the Fermi points, $\alpha$ and $\psi$ have cusps of the form
$|p\mp 1/4|^{\eta-1}$ and $|p\mp 1/4|^{\eta+\nu-1}$.
Figures~\ref{fig:lemmasolution}~(c) and (d) show the triplet solution of
Remark~\ref{rem:triplet_1d} for the same $\eta$, $\nu$ and $U_0$.
The functions $\alpha$ and $d$ are odd, $|\alpha|\le 1/2$ with
equality only at $p=\pm 1/4$,
$|d|\le 0.140$, and $\lambda=d/(2\alpha)\in[0.140, 0.267]$.

\section{Numerical Results}
\label{sec:benchmarks}

\subsection{Benchmarks against Analytical Solutions}
We use the analytical solution from Section~\ref{sec:analytical}
to study the dependency of the convergence order on the spline degree $\kappa$
and the smothness parameters $\eta$ and $\nu$ of the solution under
mesh refinement. The arising Epstein zeta functions are evaluated using
\texttt{EpsteinLib}~\cite{buchheit2024epstein}.

In Figure~\ref{fig:rates-l2-h1} we turn to a smooth reference
solution with $\eta = 4.5$
and $\nu = 1.01$. In this case, $\psi \in C^{4.51}$, so for spline
degree $\kappa$ below 5, we observe the optimal $L^2$ convergence
rate of $\kappa + 1$ in Figure~\ref{fig:rates-l2-h1}~(a), while for
$\kappa = 5$, the rate is limited by the smoothness of the solution.
With the $H^1$ norm in Figure~\ref{fig:rates-l2-h1}~(b), we lose on
order of convergence compared with the $L^2$ results. This is again
in agreement with the approximation properties of the B-splines.
Finally, in Figure~\ref{fig:rates-l2-h1}~(c) we fix the spline degree
$\kappa = 3$ and study the influence of the smoothness on the convergence rate.
From the explicit form in terms of Zeta functions the optimal rate
is inferred as $\eta + \nu - 1/2$.
This is the rate we see in all but one of the four lines.
Only for $\nu = 1.01$ the smoothness of the spline is limiting the order.

\pgfplotsset{
  convaxis/.style={
    width=6.4cm, height=5.6cm,
    log basis x=2, xmin=12, xmax=680,
    xtick={16,32,64,128,256,512}, xticklabels={16,32,64,128,256,512},
    ymin=1e-10, ymax=3e-2,
    xlabel={$n$}, ylabel={relative $L^2$ error},
    label style={font=\footnotesize},
    tick label style={font=\scriptsize},
    legend style={font=\scriptsize, draw=none, fill=none},
    legend cell align=left, legend pos=south west,
  },
  p1/.style={black, mark=o, mark size=1.5pt},
  p3/.style={blue!70!black, mark=square, mark size=1.5pt},
  p5/.style={red!70!black, mark=triangle, mark size=1.8pt},
  p7/.style={green!50!black, mark=diamond, mark size=1.8pt},
  nuA/.style={orange!80!black, mark=o, mark size=1.5pt},
  nuB/.style={violet, mark=square, mark size=1.5pt},
  nuC/.style={teal, mark=triangle, mark size=1.8pt},
  nuD/.style={magenta!80!black, mark=diamond, mark size=1.8pt},
}
\tikzset{
  slope/.style={gray, thin},
  slopelabel/.style={gray, font=\scriptsize},
}
\newcommand{\ConvergenceS}{%
\begin{tikzpicture}
  \begin{axis}[xmode=log, ymode=log, convaxis]
    \addplot[p1] coordinates {(16,4.101e-3) (32,8.025e-4) (64,1.917e-4)
      (128,4.742e-5) (256,1.183e-5) (512,2.955e-6)};
    \addlegendentry{$\kappa=1$}
    \addplot[p3] coordinates {(16,1.030e-3) (32,6.516e-5) (64,2.014e-6)
      (128,7.607e-8) (256,4.218e-9)};
    \addlegendentry{$\kappa=3$}
    \addplot[p5] coordinates {(16,9.469e-4) (32,4.718e-5) (64,6.884e-7)
      (128,2.300e-9)};
    \addlegendentry{$\kappa=5$}
    \addplot[p7] coordinates {(16,9.461e-4) (32,4.321e-5) (64,5.313e-7)
      (128,5.477e-10)};
    \addlegendentry{$\kappa=7$}
    \draw[slope] (axis cs:128,6e-6) -- (axis cs:256,1.5e-6) -- (axis cs:128,1.5e-6) -- cycle;
    \node[slopelabel, anchor=east] at (axis cs:128,3e-6) {$2$};
    \draw[slope] (axis cs:256,3e-9) -- (axis cs:512,1.875e-10) -- (axis cs:256,1.875e-10) -- cycle;
    \node[slopelabel, anchor=east] at (axis cs:256,7.5e-10) {$4$};
  \end{axis}
\end{tikzpicture}}
\newcommand{\ConvergenceD}{%
\begin{tikzpicture}
  \begin{axis}[xmode=log, ymode=log, convaxis]
    \addplot[p1] coordinates {(16,6.130e-3) (32,1.491e-3) (64,3.694e-4)
      (128,9.235e-5) (256,2.298e-5) (512,5.745e-6)};
    \addlegendentry{$\kappa=1$}
    \addplot[p3] coordinates {(16,3.876e-4) (32,3.410e-5) (64,3.025e-6)
      (128,2.698e-7) (256,2.421e-8) (512,2.203e-9)};
    \addlegendentry{$\kappa=3$}
    \addplot[p5] coordinates {(16,3.061e-4) (32,2.858e-5) (64,2.531e-6)
      (128,2.232e-7) (256,1.978e-8) (512,1.738e-9)};
    \addlegendentry{$\kappa=5$}
    \draw[slope] (axis cs:128,1e-5) -- (axis cs:256,2.5e-6) -- (axis cs:128,2.5e-6) -- cycle;
    \node[slopelabel, anchor=east] at (axis cs:128,5e-6) {$2$};
    \draw[slope] (axis cs:128,8e-9) -- (axis cs:256,7.07e-10) -- (axis cs:128,7.07e-10) -- cycle;
    \node[slopelabel, anchor=east] at (axis cs:128,2.4e-9) {$3.5$};
  \end{axis}
\end{tikzpicture}}
\newcommand{\ConvergenceLemma}{%
\begin{tikzpicture}
  \begin{axis}[xmode=log, ymode=log, convaxis, ymin=1e-13, xmax=2720,
    xtick={16,64,256,1024}, xticklabels={16,64,256,1024}]
    \addplot[p1] coordinates {(16,8.853e-4) (32,2.097e-4) (64,5.173e-5)
      (128,1.293e-5) (256,3.277e-6) (512,8.191e-7) (1024,2.048e-7)
      (2048,5.119e-8)};
    \addlegendentry{$\kappa=1$}
    \addplot[p3] coordinates {(16,1.795e-5) (32,9.440e-7) (64,5.469e-8)
      (128,3.331e-9) (256,2.066e-10) (512,1.288e-11)};
    \addlegendentry{$\kappa=3$}
    \addplot[p5] coordinates {(16,1.491e-6) (32,4.789e-8) (64,1.485e-9)
      (128,4.604e-11) (256,1.429e-12)};
    \addlegendentry{$\kappa=5$}
    \draw[slope] (axis cs:512,3e-7) -- (axis cs:1024,7.5e-8) -- (axis cs:512,7.5e-8) -- cycle;
    \node[slopelabel, anchor=east] at (axis cs:512,1.5e-7) {$2$};
    \draw[slope] (axis cs:64,4e-6) -- (axis cs:128,2.5e-7) -- (axis cs:128,4e-6) -- cycle;
    \node[slopelabel, anchor=west] at (axis cs:128,1e-6) {$4$};
    \draw[slope] (axis cs:128,1.5e-11) -- (axis cs:256,4.69e-13) -- (axis cs:128,4.69e-13) -- cycle;
    \node[slopelabel, anchor=east] at (axis cs:128,2.6e-12) {$5$};
  \end{axis}
\end{tikzpicture}}
\newcommand{\ConvergenceLemmaHone}{%
\begin{tikzpicture}
  \begin{axis}[xmode=log, ymode=log, convaxis, ymin=1e-11, ymax=3e-1,
    xmax=2720, xtick={16,64,256,1024}, xticklabels={16,64,256,1024},
    ylabel={relative $H^1$ error}]
    \addplot[p1] coordinates {(16,9.436e-2) (32,4.670e-2) (64,2.329e-2)
      (128,1.163e-2) (256,5.816e-3) (512,2.908e-3) (1024,1.454e-3)
      (2048,7.269e-4)};
    \addlegendentry{$\kappa=1$}
    \addplot[p3] coordinates {(16,1.419e-3) (32,1.634e-4) (64,1.976e-5)
      (128,2.444e-6) (256,3.046e-7) (512,3.804e-8)};
    \addlegendentry{$\kappa=3$}
    \addplot[p5] coordinates {(16,9.883e-5) (32,6.168e-6) (64,3.814e-7)
      (128,2.364e-8) (256,1.467e-9) (512,9.138e-11)};
    \addlegendentry{$\kappa=5$}
    \draw[slope] (axis cs:512,1e-3) -- (axis cs:1024,5e-4) -- (axis cs:512,5e-4) -- cycle;
    \node[slopelabel, anchor=east] at (axis cs:512,7.1e-4) {$1$};
    \draw[slope] (axis cs:64,1.6e-4) -- (axis cs:128,2e-5) -- (axis cs:128,1.6e-4) -- cycle;
    \node[slopelabel, anchor=west] at (axis cs:128,5.7e-5) {$3$};
    \draw[slope] (axis cs:128,8e-9) -- (axis cs:256,5e-10) -- (axis cs:128,5e-10) -- cycle;
    \node[slopelabel, anchor=east] at (axis cs:128,2e-9) {$4$};
  \end{axis}
\end{tikzpicture}}
\newcommand{\ConvergenceLemmaNu}{%
\begin{tikzpicture}
  \begin{axis}[xmode=log, ymode=log, convaxis, ymin=1e-14, xmax=2720,
    xtick={16,64,256,1024}, xticklabels={16,64,256,1024}]
    \addplot[nuA] coordinates {(16,4.780e-4) (32,1.010e-4) (64,2.122e-5)
      (128,4.461e-6) (256,9.389e-7) (512,1.974e-7) (1024,4.149e-8)
      (2048,8.723e-9)};
    \addlegendentry{$\nu=0.25$, $n^{-2.25}$}
    \addplot[nuB] coordinates {(16,1.681e-4) (32,2.780e-5) (64,4.581e-6)
      (128,7.555e-7) (256,1.246e-7) (512,2.056e-8) (1024,3.391e-9)
      (2048,5.593e-10)};
    \addlegendentry{$\nu=0.6$, $n^{-2.6}$}
    \addplot[nuC] coordinates {(16,4.230e-5) (32,5.182e-6) (64,6.396e-7)
      (128,7.939e-8) (256,9.896e-9) (512,1.240e-9) (1024,1.569e-10)
      (2048,2.021e-11)};
    \addlegendentry{$\nu=1.01$, $n^{-3.01}$}
    \addplot[nuD] coordinates {(16,1.840e-3) (32,1.134e-3) (64,6.978e-4)
      (128,4.293e-4) (256,2.642e-4) (512,1.626e-4) (1024,1.001e-4)
      (2048,6.161e-5)};
    \addlegendentry{$\eta=1.1$, $\nu=0.1$, $n^{-0.7}$}
  \end{axis}
\end{tikzpicture}}

\pgfplotsset{
  zoomaxis/.style={
    width=6.1cm, height=7cm,
    xmin=0.23, xmax=0.27,
    xtick={0.23,0.24,0.25,0.26,0.27},
    xticklabels={$0.23$,$0.24$,$\frac14$,$0.26$,$0.27$},
    xlabel={$p$},
    label style={font=\footnotesize},
    tick label style={font=\scriptsize},
    legend style={font=\scriptsize, draw=none, fill=none},
    legend cell align=left, legend pos=north west,
    yticklabel style={/pgf/number format/fixed, /pgf/number format/precision=2},
  },
  zoominset/.style={
    width=3.1cm, height=3.0cm,
    xmin=1e-6, xmax=0.25, ymin=1e-12, ymax=1,
    xtick={1e-5,1e-2}, ytick={1e-12,1e-8,1e-4,1},
    xlabel={$\delta$}, xlabel style={yshift=4pt},
    tick label style={font=\tiny}, label style={font=\scriptsize},
    axis background/.style={fill=white},
  },
  sref/.style={blue!70!black, thick},
  sn16/.style={blue!35, dotted, thick},
  sn64/.style={blue!55, dashed},
  sn256/.style={blue!80!black, dashdotted},
  pref/.style={red!70!black, thick},
  pn16/.style={red!35, dotted, thick},
  pn64/.style={red!55, dashed},
  pn256/.style={red!80!black, dashdotted},
}
\newcommand{\LemmaZoomPanel}[4]{%
\begin{tikzpicture}
  \begin{axis}[zoomaxis, ymin=#2, ymax=#3, ylabel={#4}]
    \draw[gray, dotted] (axis cs:0.25,#2) -- (axis cs:0.25,#3);
    \addplot[#1ref] table[x=p, y=ref] {Figures/zoom-main-#1.dat};
    \addlegendentry{exact}
    \addplot[#1n16] table[x=p, y=n0] {Figures/zoom-main-#1.dat};
    \addlegendentry{$n=16$}
    \addplot[#1n64] table[x=p, y=n1] {Figures/zoom-main-#1.dat};
    \addlegendentry{$n=64$}
    \addplot[#1n256] table[x=p, y=n2] {Figures/zoom-main-#1.dat};
    \addlegendentry{$n=256$}
    \coordinate (inset) at (axis cs:0.27,#3);
  \end{axis}
  \begin{axis}[xmode=log, ymode=log, zoominset, at={(inset)},
    anchor=outer north east, xshift=-2pt, yshift=-2pt]
    \addplot[#1ref] table[x=d, y=ref] {Figures/zoom-log-#1.dat};
    \addplot[#1n16] table[x=d, y=n0] {Figures/zoom-log-#1.dat};
    \addplot[#1n64] table[x=d, y=n1] {Figures/zoom-log-#1.dat};
    \addplot[#1n256] table[x=d, y=n2] {Figures/zoom-log-#1.dat};
  \end{axis}
\end{tikzpicture}}
\newcommand{\LemmaZoomS}{\LemmaZoomPanel{s}{0.7}{0.76}{$\psi$}}
\newcommand{\LemmaZoomP}{\LemmaZoomPanel{p}{0.04}{0.2}{$d$}}

\begin{figure}[ht]
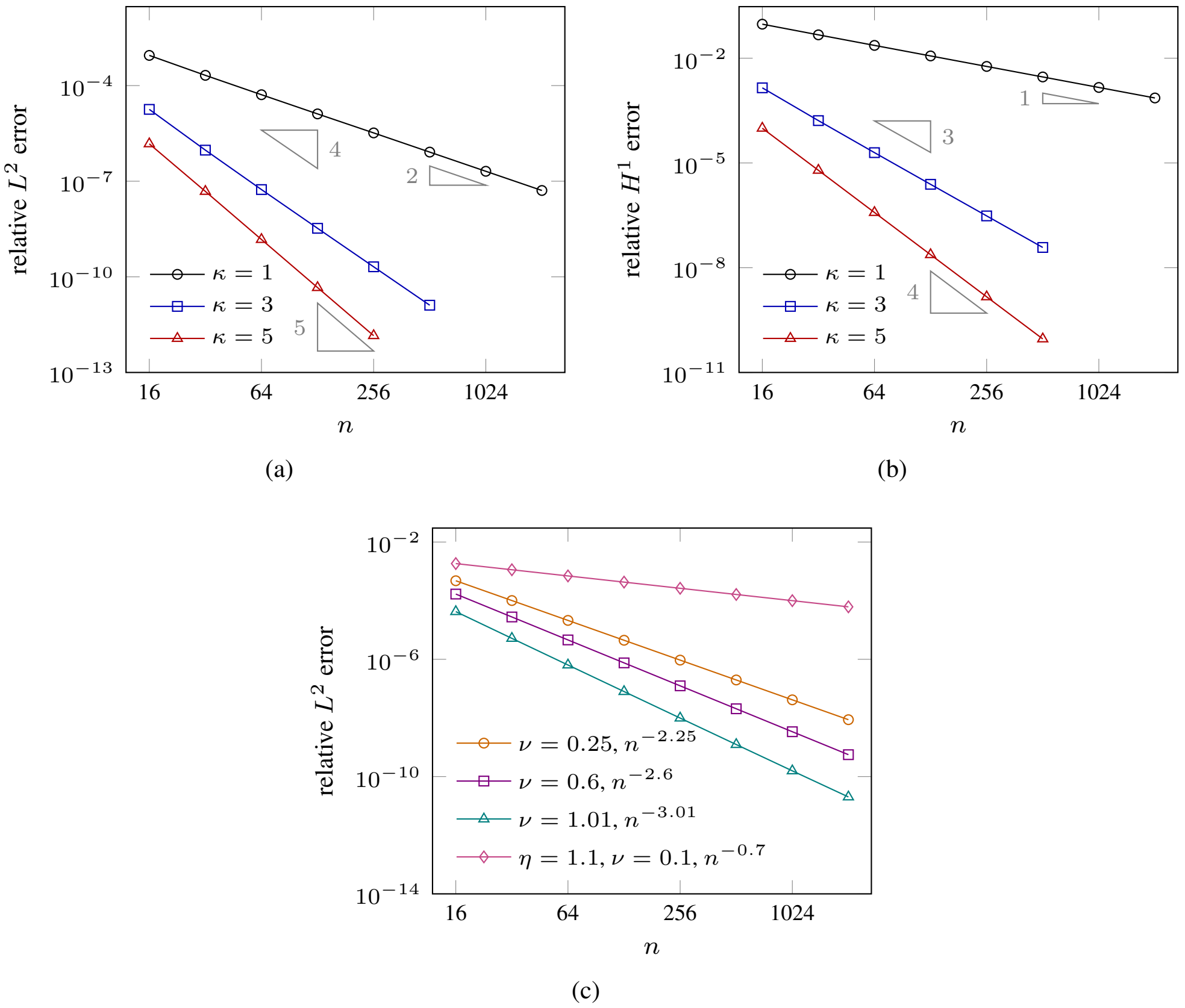

  \centering
  \subfloat[]{\ConvergenceLemma}\hfill
  \subfloat[]{\ConvergenceLemmaHone}\\
  \subfloat[]{\ConvergenceLemmaNu}
  \caption{Relative error of $\psi$ for the solution of
    Lemma~\ref{lem:singlet1d} with $C_0=2$, $U_0=0.5$, $\alpha_0=0.35$
    against the number $n$ of splines. (a) and (b): $\eta=4.5$,
    $\nu=1.01$. Gray triangles: slopes
    $n^{-s}$. (a) $L^2$ error. For $\kappa=5$, the smoothness of the
    solution limits the order. (b) $H^1$ error; for $\kappa=5$, the
    order is $\eta+\nu-3/2=4.01$ instead of $5$. (c) $L^2$ error for
    $\kappa=3$, $\eta=2.5$ and three kernel exponents $\nu$, and for
    the solution of Figure~\ref{fig:lemmasolution}, $\eta=1.1$,
  $\nu=0.1$; legend: predicted rate $n^{-(\eta+\nu-1/2)}$.}
  \label{fig:rates-l2-h1}
\end{figure}

\begin{figure}[ht]
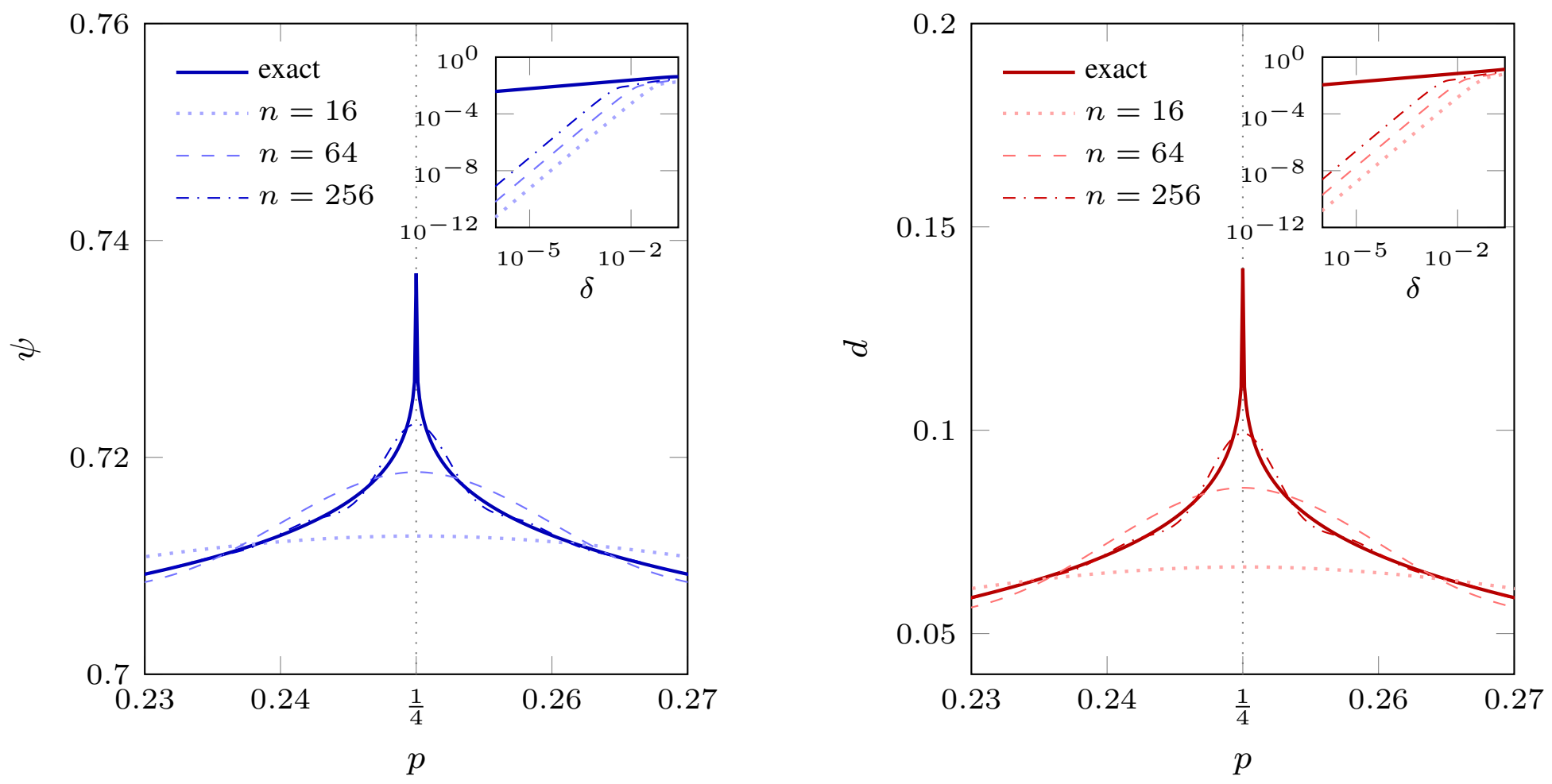

  \centering
  \subfloat[]{\LemmaZoomS}\hfill
  \subfloat[]{\LemmaZoomP}
  \caption{B-spline solutions with $\kappa=3$ near the Fermi point
    $p=1/4$ for the solutions of Figure~\ref{fig:lemmasolution}.
    (a) singlet gap $f=\psi$, (b) triplet amplitude $f=d$. Inset in both:
    $f(1/4)-f(1/4+\delta)$ against $\delta$ for the exact solution
    and for the splines on a double logarithmic scale; the exact
    solution has the cusp
    $\delta^{\eta+\nu-1}=\delta^{0.2}$, the splines are smooth at
  $p=1/4$ and decay like $\delta^2$ for $\delta<1/n$.}
  \label{fig:lemmazoom}
\end{figure}

Figure~\ref{fig:lemmazoom} shows the B-spline solution for $\kappa =
3$ near the Fermi point $p=1/4$ together with the analytical result.
The spline solutions are smooth at $p=1/4$, so that the quantity $\psi_n(1/4)-\psi_n(1/4+\delta)$
decays like $\delta^2$ for $\delta<1/n$, while the exact solution has the
cusp $\delta^{0.2}$.

\subsection{Optimal Scaling of Numerical Work}

Figure~\ref{fig:scaling} shows the work per iteration and the peak
memory of our implementation against the number $N=n^d$ of grid
cells, for cubic B-splines on the square lattice ($d=2$) and on the
fcc lattice of Section~\ref{sec:experiments} ($d=3$).
The numerical work clearly follows a linear scaling with $N$.
A linear fit over the data range gives a slope between $0.99$ and $1.01$,
both with single and multi core performance.
On the largest 3D grid, $n=256$ ($N\approx 1.7\cdot 10^7$), one
iteration with $\kappa=3$ takes $224$\,s on one core and $16.5$\,s on
$16$ cores, a speed-up of $13.6$.
The reported rimings are medians of three runs on an AMD Ryzen~9
9950X3D with $16$ cores.

\pgfplotsset{
  scalingaxis/.style={
    width=6.4cm, height=5.6cm,
    xmin=150, xmax=3e7,
    xlabel={$N=n^d$},
    label style={font=\footnotesize},
    tick label style={font=\scriptsize},
    legend style={font=\scriptsize, draw=none, fill=none},
    legend cell align=left,
  },
  d2/.style={blue!70!black, mark=o, mark size=1.5pt},
  d3/.style={red!70!black, mark=square, mark size=1.5pt},
  ref/.style={black, dashed, thin, domain=150:3e7, samples=40},
  deg3/.style={blue!70!black, mark size=1.5pt},
  deg5/.style={red!70!black, mark size=1.5pt},
  deg7/.style={green!50!black, mark size=1.5pt},
}
\newcommand{\ScalingTime}{%
\begin{tikzpicture}
  \begin{axis}[xmode=log, ymode=log, scalingaxis, ylabel={time per iteration [s]},
      ymax=1e3, legend pos=north west]
    \addplot[d2] coordinates {(256,0.54e-3) (1024,2.08e-3) (4096,8.27e-3)
      (16384,27.33e-3) (65536,111.12e-3) (262144,444.88e-3) (1048576,1.79915)
      (4194304,8.00882) (16777216,34.20174)};
    \addlegendentry{$d=2$}
    \addplot[d2, densely dotted, mark options={solid}, forget plot] coordinates {(256,0.19e-3)
      (1024,0.30e-3) (4096,0.81e-3) (16384,2.73e-3) (65536,10.40e-3)
      (262144,35.72e-3) (1048576,148.99e-3) (4194304,758.79e-3)
      (16777216,3.39376)};
    \addplot[d3] coordinates {(4096,53.85e-3) (13824,182.31e-3)
      (32768,429.31e-3) (110592,1.44579) (262144,3.47135) (884736,11.63647)
      (2097152,27.75317) (7077888,93.48032) (16777216,224.23984)};
    \addlegendentry{$d=3$}
    \addplot[d3, densely dotted, mark options={solid}, forget plot] coordinates {(4096,4.72e-3)
      (13824,13.85e-3) (32768,30.64e-3) (110592,101.59e-3) (262144,242.77e-3)
      (884736,816.60e-3) (2097152,1.97890) (7077888,6.79309)
      (16777216,16.49608)};
    \addplot[ref] {5.569e-7*x*ln(x)/ln(2)};
    \addlegendentry{$\propto N\log N$}
  \end{axis}
\end{tikzpicture}}
\newcommand{\ScalingDegree}{%
\begin{tikzpicture}
  \begin{axis}[xmode=log, ymode=log, scalingaxis,
      ylabel={time per iteration and cell [$\mu$s]},
      ymin=1, ymax=300, legend pos=north west]
    \addplot[deg3, mark=o] coordinates {(256,2.109) (1024,2.031) (4096,2.019)
      (16384,1.668) (65536,1.696) (262144,1.697) (1048576,1.716)
      (4194304,1.910) (16777216,2.039)};
    \addlegendentry{$\kappa=3$}
    \addplot[deg5, mark=o] coordinates {(256,5.195) (1024,5.117) (4096,5.068)
      (16384,4.338) (65536,4.364) (262144,4.411) (1048576,4.713)
      (4194304,5.191) (16777216,5.332)};
    \addlegendentry{$\kappa=5$}
    \addplot[deg7, mark=o] coordinates {(256,10.94) (1024,10.76) (4096,9.194)
      (16384,9.196) (65536,9.197) (262144,9.203) (1048576,9.987)
      (4194304,10.17)};
    \addlegendentry{$\kappa=7$}
    \addplot[deg3, mark=square, forget plot] coordinates {(4096,13.15)
      (13824,13.19) (32768,13.10) (110592,13.07) (262144,13.24) (884736,13.15)
      (2097152,13.23) (7077888,13.21) (16777216,13.37)};
    \addplot[deg5, mark=square, forget plot] coordinates {(4096,16.94)
      (13824,17.31) (32768,17.26) (110592,17.17) (262144,17.17) (884736,18.44)
      (2097152,21.67)};
    \addplot[deg7, mark=square, forget plot] coordinates {(4096,22.54)
      (13824,22.73) (32768,22.63) (110592,22.55) (262144,22.60) (884736,27.75)
      (2097152,32.21)};
  \end{axis}
\end{tikzpicture}}
\newcommand{\ScalingMemory}{%
\begin{tikzpicture}
  \begin{axis}[xmode=log, ymode=log, scalingaxis, ylabel={peak memory [MiB]},
      ymin=100, ymax=1e5, legend pos=north west]
    \addplot[d2] coordinates {(256,281) (1024,283) (4096,291) (16384,323)
      (65536,429) (262144,858) (1048576,2556) (4194304,9116) (16777216,35534)};
    \addlegendentry{$d=2$}
    \addplot[d3] coordinates {(4096,310) (13824,339) (32768,388) (110592,570)
      (262144,917) (884736,2336) (2097152,4900) (7077888,15585)
      (16777216,36255)};
    \addlegendentry{$d=3$}
    \addplot[ref] {349+2.143e-3*x};
    \addlegendentry{$r_0+r_1N$}
  \end{axis}
\end{tikzpicture}}

\begin{figure}[ht]
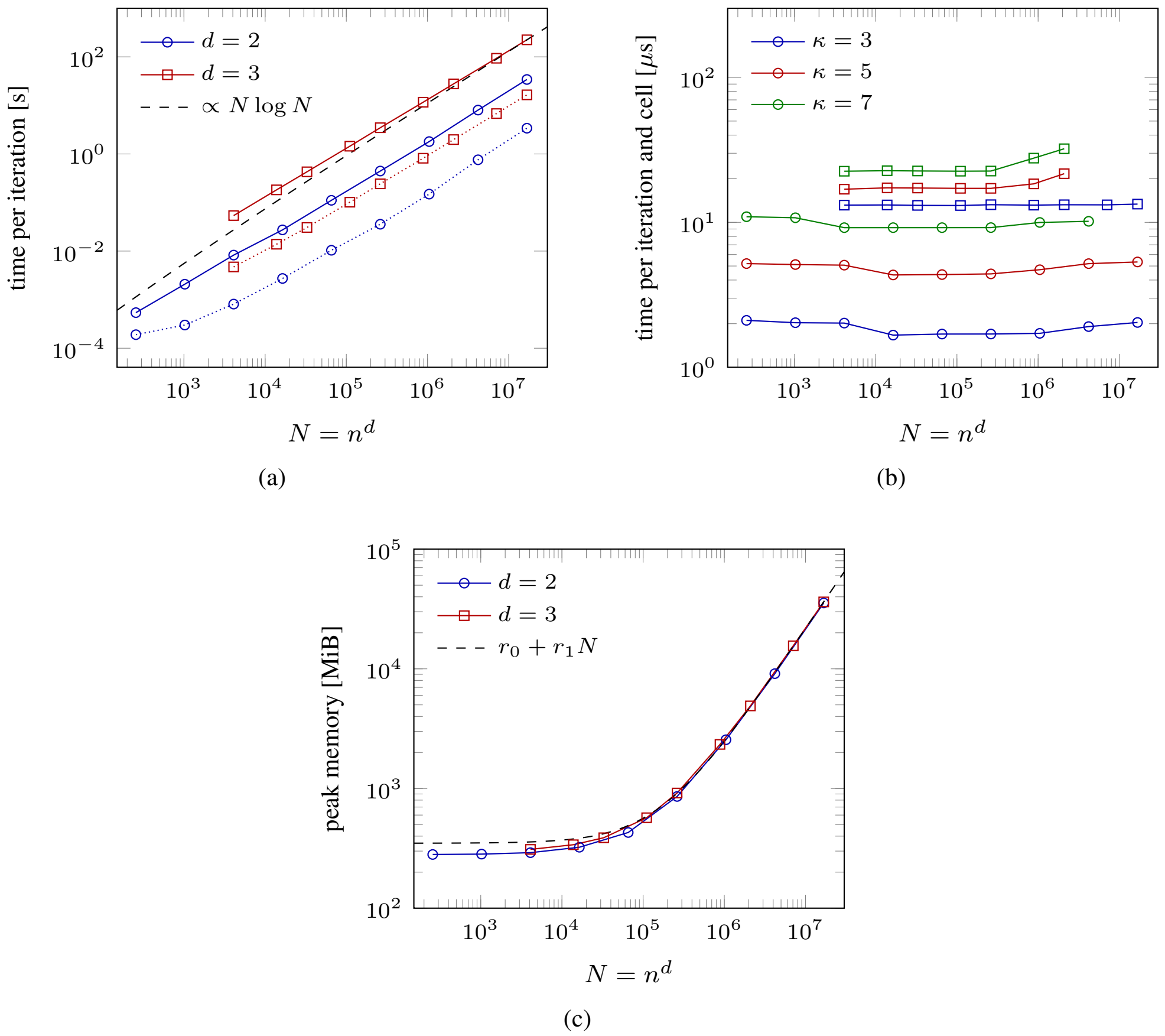

  \centering
  \subfloat[]{\ScalingTime}\hfill
  \subfloat[]{\ScalingDegree}\\
  \subfloat[]{\ScalingMemory}
  \caption{Work and memory against the number $N=n^d$ of grid cells.
    (a) Time per quasi-Newton iteration for $\kappa=3$ on $1$ core
    (solid) and $16$ cores (dotted). (b) Time per iteration and cell on
    $1$ core for $\kappa=3,5,7$; circles: $d=2$, squares: $d=3$. (c)
  Peak memory for $\kappa=3$. Dashed lines: $N\log N$ and $N$.}
  \label{fig:scaling}
\end{figure}

\subsection{Convergence Studies on Two-Dimensional Lattices}

\begin{figure}[ht]
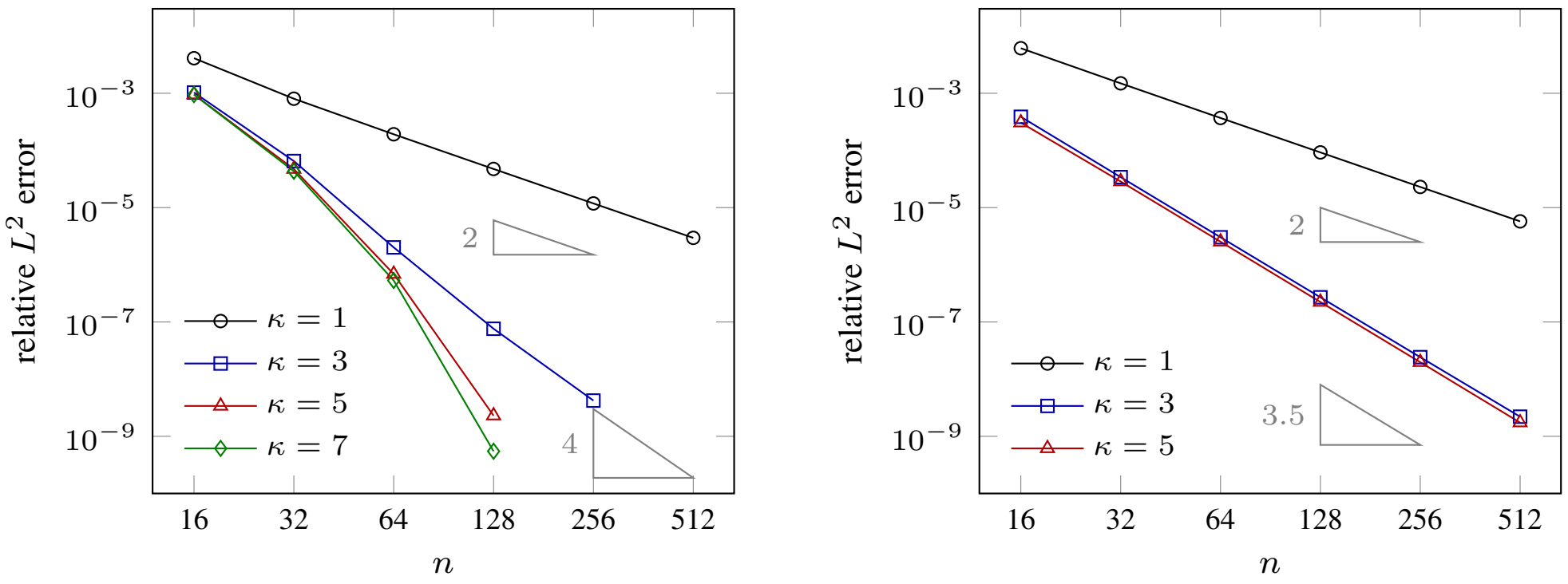

  \centering
  \subfloat[]{\ConvergenceS}\hfill
  \subfloat[]{\ConvergenceD}
  \caption{Relative $L^2$ error of $\Delta$ against the number $n$ of
    splines per dimension, $T = 0$, $\mu = 0$. Gray
    triangles: slopes $n^{-s}$. (a) $s$-wave state, triangular lattice,
    $\nu=2$, $U_0=1$, $C_0=1.5$. (b) $d$-wave,
  square lattice, $\nu=2.5$, $U_0=0.7$, $C_0=0.75$.}
  \label{fig:convergence}
\end{figure}

Here, we focus on the convergence of the B-spline methods for a fully
gapped $s$-wave solution on a triangular lattice and a nodal superconducter
with $d$-wave symmetry on a square lattice. In case of the $s$-wave,
the reference was obtained by a spectral Fourier method with $2048^2$ nodes.
All tested degrees in Figure~\ref{fig:convergence}\,(a) converge at the
optimal rate. The convergence of $\kappa=7$ is limited by the
tolerance of the nonlinear solver.

While the $d$-wave solution in Figure~\ref{fig:convergence} converges
optimally for the lowest degree $\kappa = 1$, it is limited by the smoothness
of the $d$-wave solution on the Fermi surface. The error only decays
as $\nu + 1$.

\section{Numerical Simulation of Long-Range Interacting 3D Weyl Superconductors}
\label{sec:experiments}
In this section, we demonstrate our method on a 3D face-centered-cubic lattice
with complicated node structures. Lattice vectors read
\begin{equation}
  \label{eq:fcc}
  \mathbf a_1= \tfrac12\,(0,1,1)^{T},
  \mathbf a_2 = \tfrac12\,(1,0,1)^{T},\quad
  \mathbf a_3 = \tfrac12\,(1,1,0)^{T}
\end{equation}
which form the columns of the lattice matrix
$A=(\bm a_1,\bm a_2,\bm a_3)$ with $\Lambda = A\mathds Z^3$. In this
lattice, each site has 12 nearest neighbors at distance $1/\sqrt{2}$,
normalised such that the cubic lattice constant equals $1$.
The nearest-neighbor hopping dispersion $\xi$ across the twelve neighbors is
\begin{equation}
  \label{eq:dispersion}
  \xi(\bm p) = 2\big(\cos(\pi p_1)\cos(\pi p_2) + \cos(\pi
  p_2)\cos(\pi p_3) + \cos(\pi p_3)\cos(\pi p_1)\big) - \mu
\end{equation}
with chemical potential $\mu=-0.3$.
As the interaction kernel, we use a purely attractive power law
\[
  K(\bm x) = U_0 \mathcal K_\nu(\bm x),.
\]
Its lattice Fourier transform is hence given by a
single Epstein zeta function
\begin{equation}
  \label{eq:interaction_fcc}
  \hat K(\bm p)=U_0 Z_{\Lambda,\nu}(\bm p),
\end{equation}
exhibiting a singularity at $\bm p=\bm 0$ proportional
to $\log\vert\bm p\vert$.

We solve the gap equation, choosing $\psi_0(\bm p)=0.2\,\sin(\pi
p_1)\big(\sin(\pi p_2)+i\,\sin(\pi p_3)\big)$ as the starting point,
using our method with $n = 128$ splines per dimension of degree
$\kappa = 3$, yielding a total of $3\cdot 128^3=6291456$ complex unknowns.
Using quasi-Newton with Broyden updates, the method
converges within $21$ iterations to residuum $1.5\cdot 10^{-13}$, in
$69$\,s on $16$ cores.

The resulting gap function is a spin singlet solution of the form
\[
  \hat \Delta(\bm p) =
  \begin{pmatrix}
    0 & \psi(\bm p) \\
    -\psi(\bm p) & 0
  \end{pmatrix},
\]
with $\psi(-\bm p)=\psi(\bm p)$ even and with symmetries
\begin{equation}
  \label{eq:weyl_symmetry}
  \psi(\bm p) \sim  \sin(\pi p_1)\big(\sin(\pi p_2) + i\,\sin(\pi p_3)\big).
\end{equation}
Gap functions are classified by the irreducible representations of
the point group of the normal state (thus the joint symmetries of the
  hopping amplitude $\tau$, the interaction $K$, and the lattice
$\Lambda$), here the cubic group $O_h$ of the $48$ rotations and
reflections of a cube, centered at the origin, that map the fcc
lattice onto itself. The numerical gap function here belongs to the
three-dimensional representation $T_{2g}$. Its lowest lattice
harmonics are the $d$-wave basis functions with symmetries like
\[
  d_{12}(\bm p)=\sin(\pi p_1)\sin (\pi p_2),\quad d_{23}(\bm
  p)=\sin(\pi p_2)\sin(\pi p_3),\quad d_{31}(\bm p)=\sin(\pi p_3) \sin(\pi p_1).
\]
We refer the reader to \cite{Serre1977,Dresselhaus2008} for an
introduction to this classification system and to
\cite{sigrist1991phenomenological} for its application in superconductivity.

\begin{figure}
  \centering
  \subfloat[]{\includegraphics[width=0.49\linewidth]{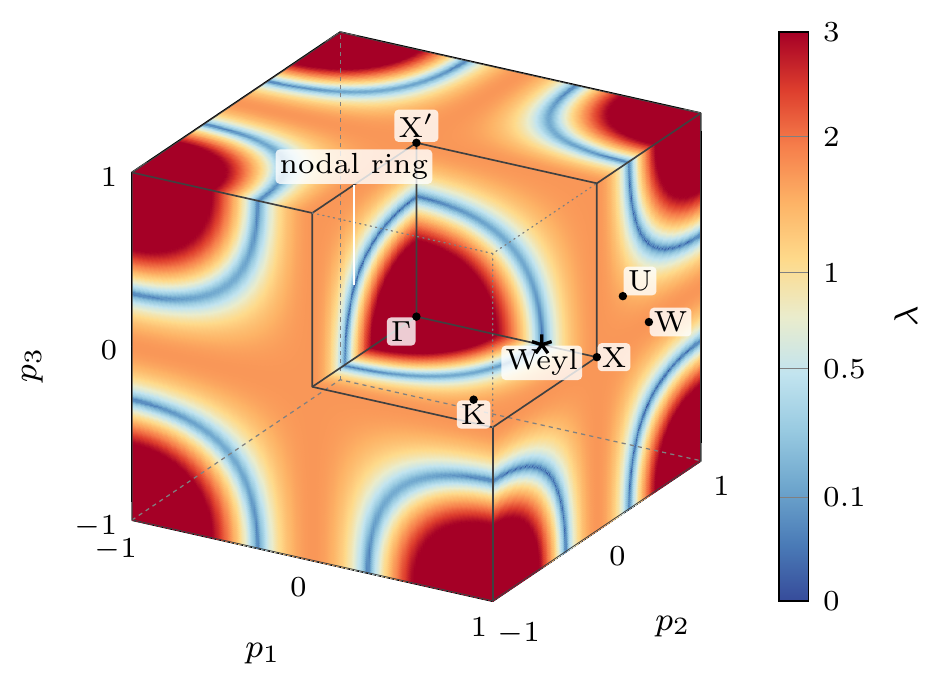}}\hfill
  \subfloat[]{\includegraphics[width=0.49\linewidth]{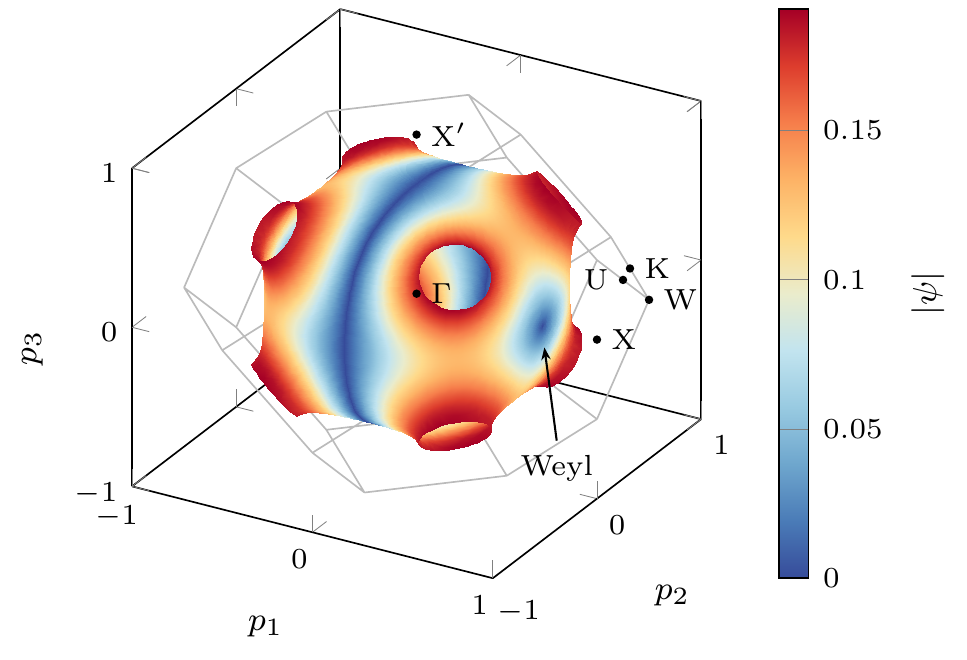}}
  \caption{Numerical solution of a chiral 3D Weyl superconductor on
    an fcc lattice with lattice vectors from Eq.~\eqref{eq:fcc} with
    nearest-neighbor hopping dispersion relation as in
    Eq.~\eqref{eq:dispersion} with chemical potential $\mu=-0.3$ and
    attractive long-range interaction as in Eq.~\eqref{eq:interaction_fcc}
    with prefactor $U_0=0.3$ and exponent $\nu=3$, $128$ cubic
    B-splines per dimension.  (a)
    Quasi-particle energies inside the BZ. The octant
    $p_1>0$, $p_2<0$, and $p_3>0$ has been exposed to better display
    the high-symmetry points and the node structure, including a ring
    node at $p_1=0$ and a Weyl node at $(\pm p_F,0,0)^\top$ with Fermi
    momentum $p_F>0$ along the $p_1$ direction. The color scale follows
    $\lambda^{1/2}$. (b) Absolute value of the singlet superconducting
    gap (upper-right element of the matrix) as a function of wave
    vector on the Fermi surface. The edges of the Brillouin zone are
  displayed as gray lines. }
  \label{fig:weyl_delta}
\end{figure}

\begin{figure}
  \centering
  \subfloat[]{\includegraphics[width=0.49\linewidth]{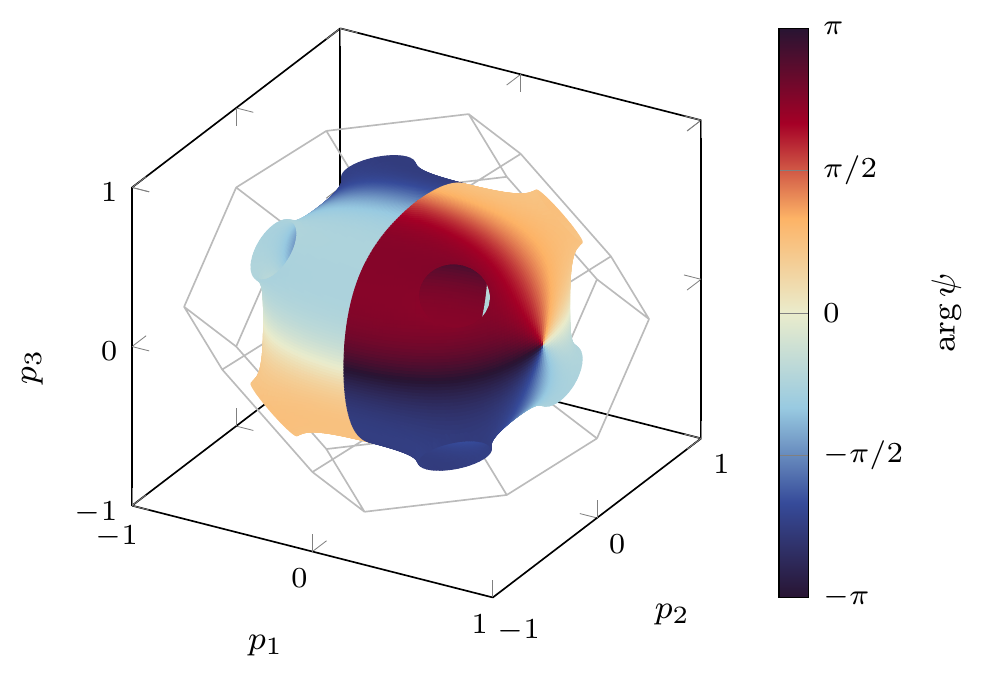}}\hfill
  \subfloat[]{\includegraphics[width=0.49\linewidth]{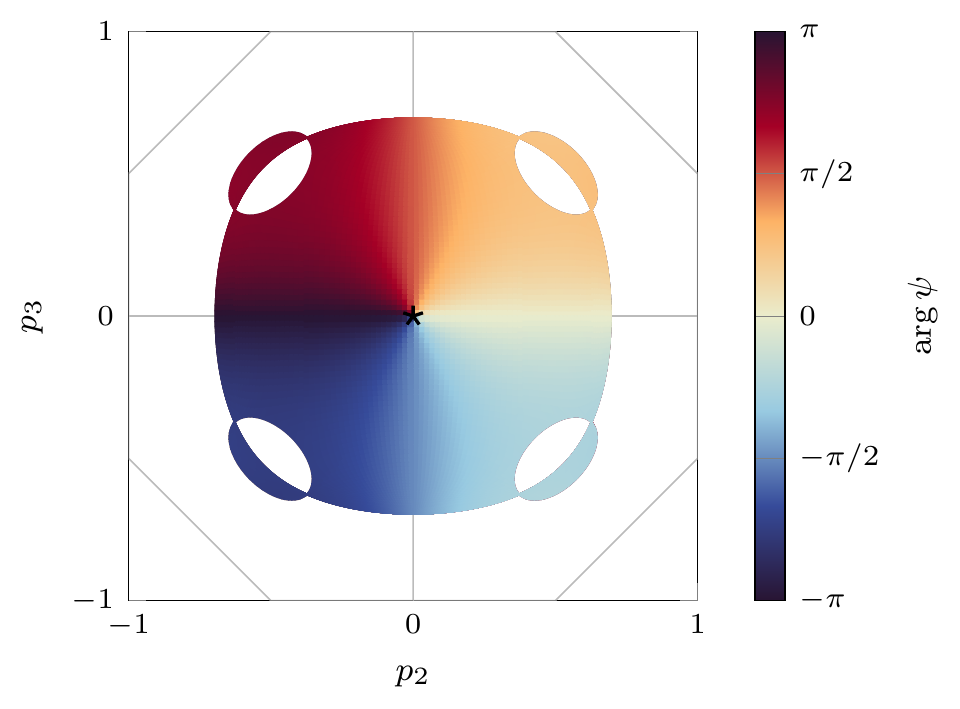}}
  \caption{Phase of singlet superconducting gap $\psi$ on the Fermi
    surface with parameters as in Fig.~\ref{fig:weyl_delta}, (b) viewed
    along the $p_1$ axis. A jump of $\pi$ occurs on the line nodes,
    whereas an integer winding around the Weyl points signals their
  topological invariant.}
  \label{fig:weyl_phase}
\end{figure}

\begin{figure}
  \centering
  \subfloat[]{\includegraphics[width=0.49\linewidth]{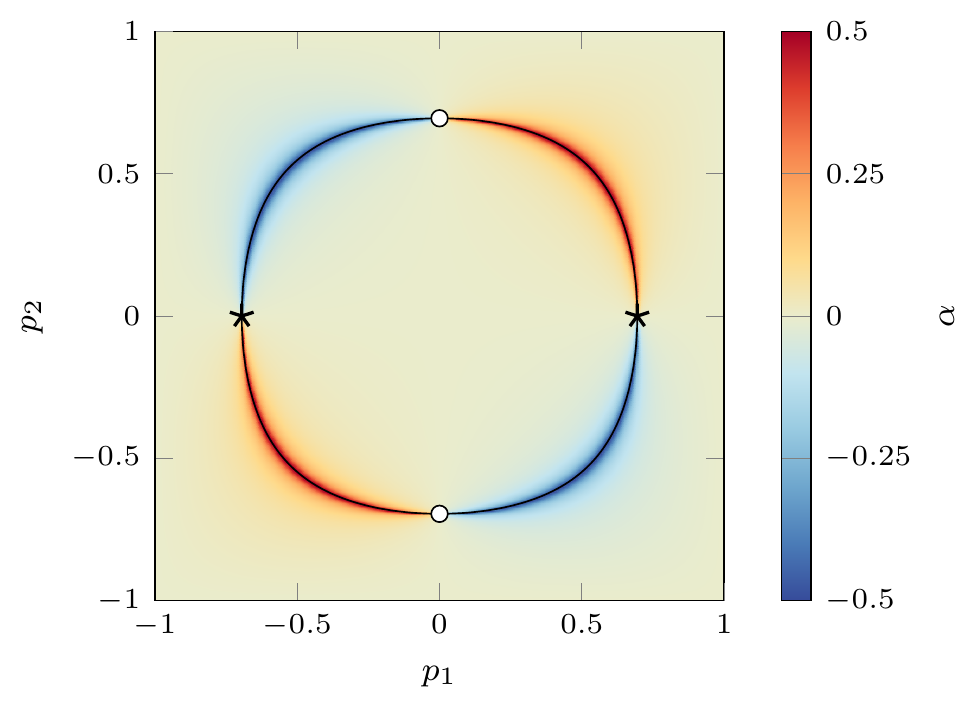}}\hfill
  \subfloat[]{\includegraphics[width=0.49\linewidth]{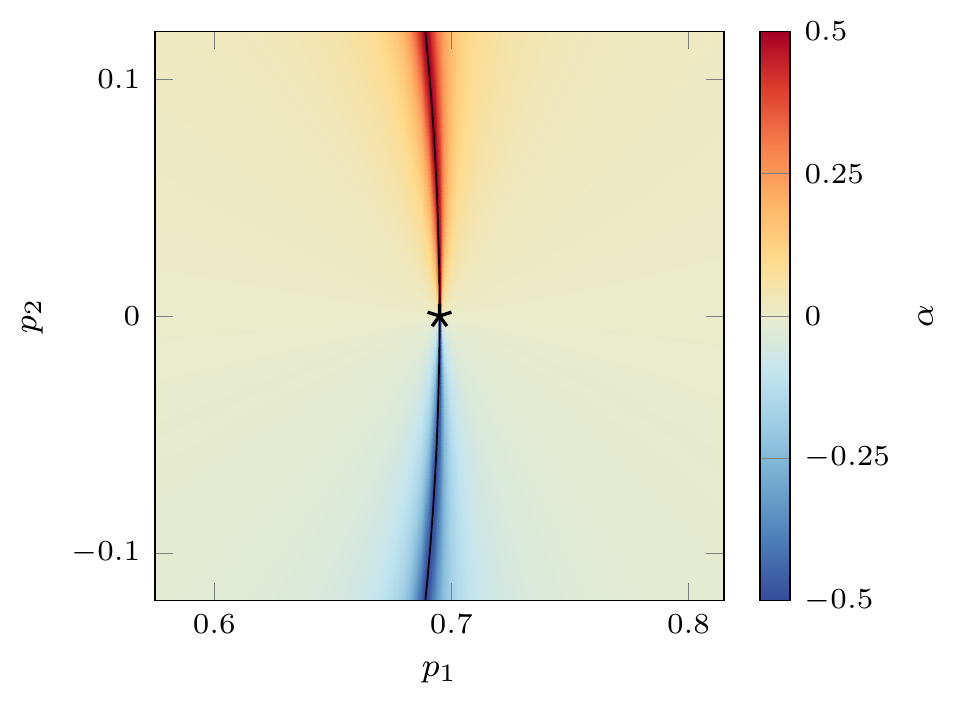}}
  \caption{Pair amplitude $\alpha=\psi/(2\lambda)$ on the plane
    $p_3=0$ (a), and near the Weyl node (b), with parameters as in
    Fig.~\ref{fig:weyl_delta}. The line marks the Fermi surface, the
    stars the Weyl nodes, and the circles the nodal ring. On the Fermi
  surface $\vert\alpha\vert=1/2$, so $\alpha$ is discontinuous at the nodes.}
  \label{fig:weyl_alpha}
\end{figure}

The structure and location of the nodes, that is wave vectors $\bm p$
at which $E(\bm p)$ is not full rank, are of particular interest.
Our numerical method allows us to resolve them faithfully.
First of all, in our case, both
eigenvalues of $E(\bm p)$ are equal,
$\lambda_1(\bm p)=\lambda_2(\bm p)=\lambda(\bm p)$,
and the quasi-particle energies obey
\[
  \lambda(\bm p) =\sqrt{\xi^2(\bm p)+\vert \psi(\bm p)\vert^2},
\]
and the node set $\mathcal N\subset E^\ast$ is the subset of the
Fermi surface $\mathcal C$ where also
$\psi(\bm p)=0$ and therefore $\lambda(\bm p)=0$ holds.
The quasi-particle energy $\lambda(\bm p)$ is displayed as a function
of $\bm p$ in Fig.~\ref{fig:weyl_delta} (a).
Meanwhile, panel (b) displays the absolute value of $\psi$ on the Fermi surface.

Two different kinds of nodes arise in our solution.
The first type is a ring node, located at the intersection of the
plane $p_1=0$ with the Fermi surface.
The second kind are the two point nodes at $(\pm p_F, 0, 0)^\top$.
The phase of $\psi(\bm p)$ winds once around the node on the Fermi
surface, with opposite sense at the two nodes. These nodes are called
Weyl nodes, as they are
similar in dispersion structure to the case of Weyl semimetals, see
\cite{schnyder2015topological,MengBalents2012}. The phase of $\psi$
on the Fermi surface as a function $\bm p$ is displayed in
Fig.~\ref{fig:weyl_phase}. Here, the integer phase vortex around the
Weyl points clearly shows their topological nature. On the other
hand, the line node displays a simple phase jump by $\pi$.

\section{Conclusions and Outlook}
\label{sec:outlook}

Unconventional superconductors are among the most well-known
materials that exhibit macroscopic quantum behavior with direct
applications in sensing, energy transport, energy storage, and
fundamental quantum research. However, the microscopic origin and
structure of the superconducting state, encoded in the
superconducting gap matrix, are unknown for many materials. Efficient
numerical simulations of unconventional superconductors are therefore
a crucial tool for investigating different interaction and pairing
models and can guide the search for exotic superconducting states. In
this work, we provide an efficient numerical method for determining
physically meaningful solutions of BCS gap equation for tight-binding
models on arbitrary Bravais lattices, a general set of attractive
interactions, including both long-range algebraic decay and compactly
supported short-ranged contributions, and both singlet and triplet
channels.  We first connect the gap equation to the physical energy functional.
We prove that the functional is bounded from below and that its
stationary states correspond to solutions of the gap equation.
Furthermore, we demonstrate that, for non-repulsive interaction
kernels, an iterative solution of the gap equation never increases
the free energy while driving the residual of the gap equation to zero.
This is an important point, as it explains why solving the gap
equation through fixed-point iteration leads to physically stable
and, therefore, meaningful results.

We combine a B-spline Galerkin discretisation of the reciprocal
lattice cell with a pseudo-differential operator approach that allows
for analytic evaluations of the arising singular convolutions,
extending the approach from \cite{buchheit2026numerical} to general
Bravais lattices and a more general set of interaction kernels. Using
the circulant structure of the arising matrices, we obtain log-linear
scaling in work with number of discretisation parameters per
iteration step and linear scaling in memory. Going further, we
determine a complete parameter family of new analytic singlet and
triplet solutions
to the gap equation for the integer lattice at zero temperature with
power-law interactions, where the gap only
decays algebraically with distance. To the best of our knowledge,
this is the first case of a non-trivial analytic solution for a BCS
gap equation with long-range interactions. In extensive benchmarks,
we first test our method against the analytic solutions as well as
highly resolved numerical data, obtaining stable convergence, with
the convergence rate determined by the regularity of the solution and
the spline order.
Due to log-linear scaling in work and linear scaling in memory, we
are subsequently able to simulate a 3D unconventional superconductor
with attractive power-law long-range interactions on an fcc cell with
a copper-like Fermi surface. We find a chiral $d$-wave solution with
both a line node and two topological Weyl nodes where the
nonlinearity in the integrand
of the gap equation becomes non-analytic. Our approach permits the
precise and efficient evaluation of the convolution of the singular
interaction kernel with the singular nonlinearity. The arising chiral $d$-wave
solution has the same symmetry as the suggested pairing for the
material URu$_2$Si$_2$ \cite{Mydosh2020}.

In future work, it will be meaningful to extend the treatment from
positive interaction kernels to kernels of general sign, employing a
more general energy minimisation approach. Furthermore,
$hp$-adaptivity can be used to improve convergence rates for
solutions with low regularity. The method can also be
generalised to multi-band and flat-band structures, where the
increase in dimension and more severe singularity structure require
care. In particular, an extension to systems with non-trivial quantum
geometry, such as twisted bilayer graphene, would be valuable
\cite{torma2022superconductivity}. Finally, an extension to the full
Hartree-Fock-Bogoliubov
theory is meaningful, as well as an extension to Eliashberg theory
\cite{Esterlis2026}. Meanwhile, combining the extended method with
realistic interaction potentials and dispersion relations could
enable quantitative predictions of superconducting material properties.

\section*{Acknowledgement}
This paper is dedicated to Professor Olaf Steinbach of Graz
University of Technology, with fond memories of decades of friendship
and collaboration with SR, beginning with B-splines.

AB and TK acknowledge  the support they received from the Klaus-Tschira Stiftung under
Grant No. 00.025.2025. The authors gratefully acknowledge the
scientific support and HPC
resources provided by the Erlangen National High Performance
Computing Center (NHR@FAU) of the Friedrich-Alexander-Universität
Erlangen-Nürnberg (FAU) under the NHR project n101af.

Claude Code (using models Opus 5.5, Fable 5.1, Fable 5.2) was used in
the preparation of this manuscript for literature research,
scientific discussion, drafting of text passages, informal proof
checking, data analysis, specification driven development, and
visualisation. All results were carefully verified by the authors.
\appendix
\section{Technical Lemmas}
In the following lemma, we prove properties of the anomalous density
$G_\beta[\hat \Delta](\bm p)$, viewed as a pointwise function of
$\hat \Delta(\bm p)$.
\begin{lemma}[Properties of anomalous density]
  \label{lem:aux_lemma_gb}
  Let $0<\beta<\infty$. For $A\in \mathds C^{2\times 2}$, define the
  function $\tilde G_\beta:\mathds C^{2\times 2}\to \mathds
  C^{2\times 2}$, by replacing in the formula for $G_\beta$ from
  Def.~\ref{def:anomalous_density} the gap $\hat \Delta$ with the
  constant matrix $A$ while fixing $\xi(\bm p)=\xi\in \mathds R$, so that
  \[
    G_\beta[\hat \Delta](\bm p)=\tilde G_\beta(\hat \Delta(\bm
    p)),\quad \xi=\xi(\bm p).
  \]
  Then $\tilde G_\beta$ is smooth and, with the derivative
  \[\delta \tilde G_\beta =  \frac{\partial}{\partial
  \varepsilon}\tilde G_\beta(A+\varepsilon B)\Big\vert_{\varepsilon=0},\]
  in direction $B$, the following properties hold for all $A,B\in
  \mathds C^{2\times 2}$.
  \begin{enumerate}
    \item{(Boundedness and symmetry)}
      \[
        \Vert \tilde G_\beta(A)\Vert_F\le \tfrac{1}{\sqrt{2}},\quad
        \tilde G_\beta(A)^\top=-\tilde G_\beta(-A^\top),
      \]
    \item{(Bounds on derivative)}

      There exists $m_A>0$, such that
      \[
        \mathrm{Re}\,\mathrm{Tr}(B^\ast \delta \tilde G_\beta )\ge
        m_A \Vert B\Vert_F^2,\quad \Vert \delta \tilde
        G_\beta\Vert_F\le \tfrac{\beta}{4}\Vert B\Vert_F.
      \]
  \end{enumerate}
  In particular, $\hat \Delta \in X$ implies $G_\beta[\hat \Delta]\in X$.
\end{lemma}
\begin{proof}
  Begin by recalling from the proof of
  Thm.~\ref{thm:gap_eq_is_stationary_point} that $\phi_\beta$ is even
  and real analytic and  $\varphi_\beta(t)=\phi_\beta(\sqrt t)$ is
  smooth on $[0,\infty)$ with
  \[
    \varphi_\beta'(t) = \frac{1}{2}\tanh(\beta \sqrt{t}/2)\frac{1}{\sqrt{t}},
  \]
  hence
  \[
    \tilde G_\beta(A)=A \varphi_\beta'(\xi^2\mathds 1_2+A^\ast A).
  \]
  Form a singular value decomposition $A=U \Sigma V^\ast$ with unitary matrices
  $U,V\in \mathds C^{2\times 2}$, and
  $\Sigma=\mathrm{diag}(\sigma_1,\sigma_2)$, with singular values
  $\sigma_1,\sigma_2\ge 0$.
  Then
  \[
    \tilde G_\beta(A) =U\gamma(\Sigma) V^\ast,\quad
    \gamma(\sigma)=\frac{\sigma\tanh(\beta e/2)}{2e},\quad
    e=\sqrt{\xi^2+\sigma^2},
  \]
  with scalar and odd $\gamma\in C^\infty(\mathds R)$. As in the
  proof of Prop.~\ref{prop:boundedness}, we readily obtain the bound
  $|\gamma|\le \tfrac{1}{2}$. The derivative of $\gamma$ reads
  \[
    \gamma'(\sigma) =
    \frac{\xi^2}{e^2}\frac{\beta}{4}\frac{\tanh(\beta e/2)}{\beta
    e/2}+\frac{\sigma^2}{e^2}\frac{\beta}{4} \big(1-\tanh^2(\beta e/2)\big).
  \]
  With $x=\beta e/2$, use that $0<\tanh(x)\le x$ in the first term
  and $0\le \tanh^2(x)<1$ in the second. As the weights $\xi^2/e^2$
  and $\sigma^2/e^2$ are nonnegative and sum to one, this yields for $e>0$
  \[
    0<\gamma'(\sigma)\le \tfrac{\beta}{4}.
  \]
  For $\xi=0$ we have $\gamma(\sigma)=\tfrac12\tanh(\beta\sigma/2)$,
  hence $\gamma'(0)=\beta/4$ also at $e=0$, and the bound holds for
  all $\sigma\in\mathds R$. Furthermore, $\tilde G_\beta$ is smooth,
  as $\varphi_\beta'$ is real analytic on a neighbourhood of
  $[0,\infty)$, which contains the spectrum of $\xi^2\mathds 1_2+A^\ast A$.

  For Statement 1, note that
  \[
    \Vert \tilde G_\beta(A)\Vert_F^2 =
    \gamma(\sigma_1)^2+\gamma(\sigma_2)^2\le \tfrac{1}{2},
  \]
  hence $\Vert \tilde G_\beta(A)\Vert_F\le \tfrac{1}{\sqrt{2}}$. To
  establish the symmetry, note that  $-A^\top=(-\bar V)\Sigma (\bar
  U)^\ast$, with $\bar U$ the entrywise complex conjugate, is an SVD and thus
  \[
    \tilde G_\beta (-A^\top)=(-\bar V)\gamma(\Sigma) (\bar
    U)^\ast=-\big(U\gamma(\Sigma)V^\ast\big)^\top=-\tilde G_\beta(A)^\top.
  \]

  For Statement 2, first define the Hermitian matrices
  \[
    M=
    \begin{pmatrix}
      0 & A\\
      A^\ast & 0
    \end{pmatrix},\quad N=
    \begin{pmatrix}
      0 & B\\
      B^\ast & 0
    \end{pmatrix},
  \]
  with $\Vert N\Vert_F^2=2 \Vert B\Vert_F^2$. Since
  $M^2=\mathrm{diag}(A A^\ast, A^\ast A)$, the four eigenvalues of
  $M$ are the two singular values of $A$ equipped with positive and
  negative signs each. In particular, the eigenvalues are contained
  in the interval
  \[
    I_A=[-\Vert A\Vert_F,\Vert A\Vert_F].
  \]
  As $\gamma(x)= x \varphi_\beta'(\xi^2+x^2)$ for $x\in \mathds R$,
  we obtain using that $M^2$ is block diagonal
  \[
    \gamma(M)=M \varphi_\beta'(\xi^2 \mathds 1_4+M^2)=
    \begin{pmatrix}
      0 & A \varphi_\beta'(\xi^2 \mathds 1_2+A^\ast A) \\
      A^\ast \varphi_\beta'(\xi^2 \mathds 1_2+A A^\ast) & 0
    \end{pmatrix}.
  \]
  The upper right block is thus $\tilde G_\beta(A)$. The same holds
  with $A+\varepsilon B$, corresponding to $M+\varepsilon N$. Thus
  $\gamma(M+\varepsilon N)$ is Hermitian and block off-diagonal for
  every $\varepsilon$ with derivative
  \[
    \delta \gamma(M) = \frac{\partial}{\partial
    \varepsilon}\gamma(M+\varepsilon N)\Big\vert_{\varepsilon=0}=
    \begin{pmatrix}
      0 & \delta \tilde G_\beta\\
      (\delta \tilde G_\beta)^\ast & 0
    \end{pmatrix}.
  \]
  Therefore,
  \[
    \mathrm{Tr}(N\delta \gamma(M))=2\mathrm{Re}\,\mathrm{Tr}(B^\ast
    \delta \tilde G_\beta),\quad \text{and}\quad \Vert \delta
    \gamma(M)\Vert_F^2=2 \Vert \delta \tilde G_\beta\Vert_F^2.
  \]
  Now, diagonalise $M$ with unitary transform $Z\in \mathds C^{4\times 4}$
  \[
    M=Z \mathrm{diag}(a_1,\dots,a_4) Z^\ast,
  \]
  and define $Y=Z^\ast N Z$. The Daleckii–Krein formula
  \cite[Thm.~3.11]{Higham2008} then allows to rewrite the derivative
  in the explicit form
  \[
    \delta \gamma(M)=Z(F\odot Y)Z^\ast,
  \]
  with $\odot$ the entrywise Hadamard product and
  \[
    F_{k\ell}=\frac{\gamma(a_k)-\gamma(a_\ell)}{a_k-a_\ell},\quad
    a_k\neq a_\ell,
  \]
  with $F_{k\ell}=\gamma'(a_k)$ for $a_k= a_\ell$. By the mean value
  theorem, $F_{k\ell}=\gamma'(\theta_{k\ell})$ for some
  $\theta_{k\ell}$ between $a_k$ and $a_\ell$.

  Now, all eigenvalues lie in $I_A$, hence by the bound on $\gamma'$,
  which is continuous and positive on the compact interval $I_A$,
  \[
    m_A\le F_{k\ell} \le \frac{\beta}{4},\quad \text{with}\quad
    m_A=\min_{x \in I_A} \gamma'(x)>0.
  \]

  Using that the trace is cyclic, that unitary transformations leave
  the Frobenius norm unchanged, and that $Y^\ast =Y$, we have
  \[
    \mathrm{Tr}\big(N\delta\gamma(M)\big)=\mathrm{Tr}\Big( Y(F\odot
    Y)\Big)=\sum_{k, \ell=1}^4 Y_{k \ell} F_{\ell k} Y_{\ell k}
    =\sum_{k,\ell=1}^4  F_{\ell k} |Y_{\ell k}|^2\ge m_A \Vert N\Vert_F^2.
  \]
  Furthermore,
  \[
    \Vert \delta \gamma(M)\Vert_F^2=\Vert F\odot Y\Vert_F^2 =
    \sum_{k,\ell=1}^4  F_{\ell k}^2 |Y_{\ell k}|^2\le
    \left(\beta/4\right)^2 \Vert N\Vert_F^2.
  \]
  After inserting the block identities, the factors 2 cancel on both
  sides, and we obtain
  \[
    \mathrm{Re}\,\mathrm{Tr}\big(B^\ast \delta \tilde G_\beta\big)\ge
    m_A \Vert B\Vert_F^2,\quad \text{and}\quad \Vert \delta \tilde
    G_\beta\Vert_F \le (\beta/4) \Vert B\Vert_F.
  \]

  Finally, consider $\hat \Delta\in X$. Then, pointwise $G_\beta[\hat
  \Delta](\bm p)=\tilde G_\beta(\hat \Delta(\bm p))$ with
  $\xi=\xi(\bm p)$, which is bounded by Statement 1 and hence lies in
  $L^2(E^\ast,\mathds C^{2\times 2})$. Furthermore, since $\xi$ is
  even, the symmetry of Statement 1 applies as well and
  \[
    G_\beta[\hat \Delta](-\bm p) =\tilde G_\beta\big(\hat \Delta(-\bm
    p)\big)= -\tilde G_\beta\big(\hat \Delta(\bm p)\big)^\top=
    -G_\beta[\hat \Delta](\bm p)^\top,
  \]
  therefore $G_\beta[\hat \Delta]\in X$, concluding the proof.
\end{proof}

Based on the results for the anomalous density in the preceding
lemma, we now prove  Prop.~\ref{prop:convergence}.

\begin{proofof}{Proposition~\ref{prop:convergence}}
  We begin with the free energy in the form of
  Lem.~\ref{lem:energy_reformulation},
  \[
    F_\beta[\hat \Delta]=\tfrac{1}{2}\mathrm{tr}(\xi\mathds
    1_2)-\tfrac{1}{2} \mathrm{tr}\,\phi_\beta(E[\hat
    \Delta])+\mathrm{Re}\langle \hat \Delta,G_\beta[\hat \Delta]
    \rangle-\tfrac{1}{2}\langle G_\beta[\hat \Delta],\hat K\ast
    G_\beta[\hat \Delta]\rangle.
  \]
  Then define $\hat H$ as the difference between two consecutive gap iterates,
  \[
    \hat H=\hat \Delta_{j+1}-\hat \Delta_j,
  \]
  and, with $0\le \varepsilon\le 1$, form the convex combination
  \[
    \hat \Delta_j+\varepsilon \hat H=(1-\varepsilon)\hat
    \Delta_j+\varepsilon \hat \Delta_{j+1}.
  \]
  Define $f(\varepsilon)=\tfrac{1}{2} \mathrm{tr}\,\phi_\beta(E[\hat
  \Delta_j+\varepsilon \hat H])$. Then, its first derivative reads by
  the proof of Thm.~\ref{thm:gap_eq_is_stationary_point},
  \[
    f'(\varepsilon)=\mathrm{Re}\,\langle \hat H,G_\beta[\hat
    \Delta_j+\varepsilon \hat H] \rangle,\quad
    f'(0)=\mathrm{Re}\,\langle \hat H, G_j \rangle,
  \]
  and, by dominated convergence, we obtain for the second derivative
  \[
    f''(\varepsilon)=\frac{1}{|E^\ast|}\int\limits_{E^\ast}
    \mathrm{Re}\,\mathrm{Tr}\Big(\hat H(\bm p)^\ast \delta \tilde
    G_\beta\Big)\,\mathrm d \bm p\ge 0
  \]
  with $\tilde G_\beta$ from Lem.~\ref{lem:aux_lemma_gb} evaluated at
  $A=\hat \Delta_j(\bm p)+\varepsilon \hat H(\bm p)$ and $B=\hat
  H(\bm p)$, where the right-hand side is non-negative by the lemma's
  second statement. Using Taylor's theorem with integral remainder, we have
  \[
    f(1)-f(0)=f'(0)+\int_{0}^1 f''(t)(1-t)\,\mathrm d t\ge f'(0).
  \]
  With this result, we can bound the kinetic energy difference for
  subsequent iterates by
  \[
    \tfrac{1}{2}\mathrm{tr}\,\phi_\beta(E[\hat \Delta_{j+1}])\ge
    \tfrac{1}{2}\mathrm{tr}\,\phi_\beta(E[\hat \Delta_{j}])+
    \mathrm{Re}\,\langle \hat \Delta_{j+1}-\hat \Delta_j, G_j \rangle.
  \]
  Inserting this into the difference of free energies of subsequent
  iterates, using $\hat \Delta_{j+1}=\hat K\ast G_j$ and the
  self-adjointness of the convolution with $\hat K$, yields after
  reordering that
  \begin{equation}
    \label{eq:prop_first_inequality}
    F_\beta[\hat \Delta_{j+1}]-F_\beta[\hat \Delta_j] \le
    -\tfrac{1}{2}\langle G_{j+1}-G_j, \hat K
    \ast (G_{j+1}-G_j)\rangle,
  \end{equation}
  which is the first inequality of the statement, which holds without
  restriction on the sign of $K$. Now, as in the proof of
  Prop.~\ref{prop:boundedness}, use $0\le K\le K_\mathrm{sup}$. Then
  for general $\hat g\in L^2(E^\ast,\mathds C^{2\times 2})$ with
  $g=\mathcal F_\Lambda^{-1}\hat g$,
  \[
    \langle \hat g,\hat K\ast \hat g\rangle= \sum_{\bm x\in \Lambda}
    K(\bm x)\Vert g(\bm x)\Vert_F^2\ge 0,
  \]
  yielding
  \[
    -\tfrac{1}{2}\langle G_{j+1}-G_j, \hat K
    \ast (G_{j+1}-G_j)\rangle \le 0
  \]
  which completes the first part of the statement. For the second, notice that
  \[
    \Vert \hat K\ast \hat g\Vert^2=\sum_{\bm x\in \Lambda} K(\bm x)^2
    \Vert g(\bm x)\Vert^2_F \le K_\mathrm{sup}\langle \hat g,\hat
    K\ast \hat g\rangle.
  \]
  Writing the residual for $j\ge 1$ in convolution form,
  \[
    R_j=\hat \Delta_j-\hat K\ast G_j=\hat K\ast (G_{j-1}-G_j),
  \]
  the above inequality provides us with
  \[
    \Vert R_j\Vert^2\le K_\mathrm{sup} \langle G_{j-1}-G_j,\hat K\ast
    (G_{j-1}-G_j)\rangle\le 2 K_{\mathrm{sup}}(F_\beta[\hat
    \Delta_{j-1}]-F_\beta[\hat \Delta_{j}]),
  \]
  where the second inequality follows from
  inequality~\eqref{eq:prop_first_inequality} at iteration step $j-1$.
  Summing over $j=1,\dots,\ell$ yields
  \[
    \sum_{j=1}^\ell \Vert R_j\Vert^2\le 2
    K_{\mathrm{sup}}(F_\beta[\hat \Delta_{0}]-F_\beta[\hat \Delta_{\ell}]).
  \]
  Finally, Prop.~\ref{prop:boundedness} bounds $F_\beta[\hat
  \Delta_\ell]\ge F_\beta[0]-\tfrac{1}{4} K_\mathrm{sup}$ uniformly
  in $\ell$, and the statement follows in the limit $\ell\to \infty$.
\end{proofof}

\bibliographystyle{siam}
\bibliography{literature_new}
\end{document}